\documentclass[11pt,a4paper,reqno]{amsart}
\usepackage[T1]{fontenc}
\usepackage{lmodern,amsmath,amssymb,amsthm,mathtools,microtype}
\usepackage[margin=1in]{geometry}
\usepackage{enumitem,xcolor}
\usepackage[hypertexnames=false,colorlinks=true,linkcolor=blue!35!black,citecolor=blue!35!black,urlcolor=blue!35!black]{hyperref}
\numberwithin{equation}{section}

\newtheorem{theorem}{Theorem}[section]
\newtheorem{proposition}[theorem]{Proposition}
\newtheorem{lemma}[theorem]{Lemma}
\newtheorem{corollary}[theorem]{Corollary}
\theoremstyle{remark}
\newcommand{\E}{\mathbb E}\newcommand{\PP}{\mathbb P}
\newcommand{\Var}{\operatorname{Var}}
\newcommand{\tr}{\operatorname{tr}}\newcommand{\rank}{\operatorname{rank}}
\newcommand{\op}{\mathrm{op}}
\newcommand{\HS}{\mathcal S_2}\newcommand{\HC}{\mathrm{HC}}

\newcommand{\eps}{\varepsilon}
\newcommand{\dd}{\,\mathrm d}
\DeclareMathOperator{\diag}{diag}
\setlist[enumerate]{leftmargin=*,itemsep=2pt,topsep=4pt}
\allowdisplaybreaks[2]
\hypersetup{pdftitle={A Two-Spectrum Deterministic Equivalent for Random Skew-Kernel Matrices},pdfauthor={Guangqian Zhao}}
\title[A two-spectrum deterministic equivalent]{A Two-Spectrum Deterministic Equivalent for Random Skew-Kernel Matrices}
\author{Guangqian Zhao}
\address{School of Mathematical Sciences, University of Science and Technology of China, Hefei 230026, China}
\email{zhaoguangqian@mail.ustc.edu.cn}
\date{}
\begin{document}

\begin{abstract}
Rank-one perturbations of the mean characteristic polynomial of a random skew-kernel matrix assign a positive spatial weight to each root. For independent centered isotropic rows, the resulting response tracks the sample almost surely against each deterministic spatial test, uniformly over positive scales and arbitrary inputs, exactly when row energies concentrate on common projections of proportional rank. Directional absolute first moments determine this condition. The determinant Jensen gap controls the spatial response in every real orthonormal basis. For identity spatial covariance, logarithmic derivatives of the mean determinant give the temporal weights of the row energies it detects. Flat kernels of one prescribed rank suffice to test the row condition. A scalar equation in the two input spectra approximates the matrix response of the roots within $4/d$ in normalized trace norm. A varying ordinary response and the same-scale overlap of independent samples recover both input laws up to reciprocal rescaling, with Lipschitz stability of the relative spatial parameter and an optimal logarithmic modulus for the compactified laws on intervals of positive variation. Brownian area gives explicit spatial half-Cauchy mixtures.
\end{abstract}
\maketitle
\section{The spatial deterministic equivalent}
Which spatial directions carry the spectrum of a random alternating form? We use rank-one perturbations of the mean characteristic polynomial to locate its mass in deterministic subspaces.

Let $X_d$ have $d$ independent centered isotropic cylindrical rows on a real separable Hilbert space $\mathcal H_d$, with $\E X_{i,d}(h)X_{i,d}(g)=\langle h,g\rangle$. For real skew-adjoint Hilbert--Schmidt $K_d$ and deterministic real symmetric $Q_d\succeq0$ of size $d$, set
\[
 A_d^X=Q_d^{1/2}(X_dK_dX_d^{\mathsf T})Q_d^{1/2},\qquad
 \mu_d^{Q,X}=d^{-1}\sum_j\delta_{s_j(A_d^X)} .
\]
The middle matrix has zero diagonal and off-diagonal entries defined by $L^2$ truncation. For isonormal rows write $A_d,\mu_d^Q$. On $I_*=[0,\infty]$ use $w_t(s)=t^2s^2/(1+t^2s^2)$ with $w_t(\infty)=1$. Total mass and the responses at positive rational $t$ determine a finite measure by Stieltjes uniqueness.

Write $q_i$ for the eigenvalues of $Q_d$, $\lambda_j$ for the singular values of $K_d$ counted once per skew pair, and $\tr_d=d^{-1}\tr$. For $D=\det(I+tA_d^X)$ and $W_d(t)=w_t(|A_d^X|)$ put
\[
 P_d(t)=\E D,\qquad N_d(t)=\E[D W_d(t)]/P_d(t),\qquad
 J_d(t)=d^{-1}(\log P_d(t)-\E\log D).
\]
For $p_d(z)=\E\det(zI-iA_d^X)$, Proposition~\ref{reb:mean-rows} gives a positive matrix weight for each distinct root via
\[
 \frac{\E\operatorname{adj}(zI-iA_d^X)}{p_d(z)}
   =\sum_{a:\,p_d(a)=0}\frac{B_a}{z-a},
 \qquad B_a\succeq0,\quad \sum_a B_a=I,\quad \tr B_a=\operatorname{mult}(a).
\]
Thus $\Xi_d=\sum_a B_a\delta_{|a|}$ has normalized trace the ordinary root law $\nu_d^Q$ and response $N_d(t)=\int w_t\,\dd\Xi_d$. The row condition uses a common temporal projection $P$, with $\sup\varnothing=0$:
\begin{align}
 \Delta_d(a)&:=\sup_{ad\le p=\rank P\le d/4}
 \frac1d\sum_{i=1}^d\E\left|\|PX_i\|^2/p-1\right|,\label{eq:pc-stat}\\
 \Delta_d(a)&\longrightarrow0\qquad(0<a<1/4).\tag{PC}\label{eq:pc}
\end{align}

Let $r_d$ be the largest even integer at most $\min(d/4,\dim\mathcal H_d/2)$. On each rank-$r_d$ plane $P$ choose an orthogonal complex structure $J_P$ extended by zero, with $J_P=0$ at rank zero. Write $J_{d,P}=J_d(1)$ for $Q_d=I$, $K_d=J_P/d$.

\begin{theorem}[Spatial deterministic equivalent]\label{thm:main}
For independent centered isotropic rows, \eqref{eq:pc}, the scalar test $\sup_{\rank P=r_d}J_{d,P}\to0$, and the following assertion for every deterministic input pair $Q_d,K_d$ are equivalent:
\begin{equation}\label{eq:mean-matrix-track}
 \sup_{t>0} d^{-1}\|\E W_d(t)-N_d(t)\|_{\mathcal S_1}\longrightarrow0.
\end{equation}
Each is also equivalent to convergence in probability in the following comparison for every pair and every deterministic sequence of real symmetric contractions $B_d$:
\begin{equation}\label{eq:marked-track}
 \sup_{t>0}d^{-1}\big|\tr B_d\{W_d(t)-N_d(t)\}\big|\longrightarrow0 .
\end{equation}
Under \eqref{eq:pc}, \eqref{eq:marked-track} holds almost surely and in $L^1$ for each chosen array and test, under any cross-dimensional coupling. For the ordinary empirical law this gives
\begin{equation}\label{eq:track}
 \sup_{t>0}\left|\int w_t\,\dd\mu_d^{Q,X}-\Theta_d(t)\right|\longrightarrow0
\end{equation}
in the same modes, with $\Theta_d(t)=d^{-1}\tr N_d(t)$ and almost-sure weak tracking of $\nu_d^Q$ on $I_*$.
\end{theorem}
A spatial projection $E$ observes $d^{-1}\tr[E\Xi_d(\cdot)]$; rank $o(d)$ tests give vanishing normalized response. For Gaussian rows, $Q_d=I$ and $\lfloor d/4\rfloor$ time pairs of size $1/d$, at least half the singular values vanish while $N_d(1)$ is scalar with positive limit by Theorem~\ref{reb:mode}. This example explains the expectation in \eqref{eq:mean-matrix-track}.

To compute the comparison, put $p(u)=u/(1+u)$, $D_Q(v)=\sum_i p(vq_i)$ and $D_K(y)=2\sum_j w_y(\lambda_j)$. For nonzero inputs set $M_d(t)=p(vQ_d)$, where $v=v_d(t)>0$ solves
\begin{equation}\label{reb:two-spectral-equation}
 w=\sum_i\frac{q_i}{1+vq_i},\qquad \theta_d(t)=D_Q(v)=D_K(tw).
\end{equation}
Take $v_d=\theta_d=0$ for zero inputs.
With $e_j$ the elementary symmetric function and $e_0=1$, define
\begin{equation}\label{eq:two-spectrum-polynomials}
 a_{k,d}=(2k)!\,e_{2k}(q)e_k(\lambda^2),\qquad
 P_d(t)=\sum_{2k\le d}a_{k,d}t^{2k},\qquad
 p_d(x)=\sum_{2k\le d}(-1)^k a_{k,d}x^{d-2k}.
\end{equation}
These are the two mean determinants by Proposition~\ref{reb:mean-rows}. Pairing opposite roots gives
\begin{equation}\label{eq:two-spectrum-response}
 \Theta_d(t)=\int w_t\,\dd\nu_d^Q=\frac{tP_d'(t)}{dP_d(t)} .
\end{equation}
\begin{theorem}[Finite-dimensional comparison from two spectra]\label{reb:mode}
The crossing in \eqref{reb:two-spectral-equation} exists and is unique for nonzero inputs.
At each $0<x<\min(\rank Q_d,\rank K_d)$ the scale at which $\theta_d(t)=x$ is
\begin{equation}\label{reb:inverse-scales}
 t(x)=\frac{D_Q^{-1}(x)\,D_K^{-1}(x)}{x}.
\end{equation}
With $\theta_d=0$ if either input vanishes, the finite-dimensional root response satisfies
\begin{equation}\label{reb:two-scalar-accuracy}
 \sup_{t>0}|d\Theta_d(t)-\theta_d(t)|\le4 .
\end{equation}
For every $t>0$ and all independent centered isotropic row laws, $N_d(t)$ equals $\operatorname{diag}(q_i\partial_{q_i}\log P_d(t))$ in a $Q_d$-eigenbasis. In any real orthonormal basis the diagonal entries $m_i$ of $N_d(t)$ and $r_i$ of $\E W_d(t)$ satisfy
$\sum_i\operatorname{kl}(r_i\Vert m_i)\le 2dJ_d(t)$, where $\operatorname{kl}$ is Bernoulli relative entropy. We have
\begin{equation}\label{eq:finite-spatial-gap}
 d^{-1}\|\E W_d(t)-N_d(t)\|_{\HS}^2\le J_d(t),\qquad
 d^{-1}\|N_d(t)-M_d(t)\|_{\mathcal S_1}\le 4/d .
\end{equation}
Thus $M_d$ may replace $N_d$ throughout Theorem~\ref{thm:main}, and $d^{-1}\|\E W_d-M_d\|_{\mathcal S_1}\le\sqrt{J_d}+4/d$. The order $d^{-1}$ is optimal uniformly; the square-root dependence of the first response difference on $J_d(t)$ is also sharp (Appendix~\ref{r21:sec-logarithms}).
\end{theorem}

For a simple root $a$, differentiation under a shift $huu^*$ gives speed $u^*B_au$, the eigenvector weight in the deterministic case \cite{DPTZ}. Yaskov \cite{Y16} characterized the Marchenko--Pastur law under projections by quadratic-form concentration. Here both input spectra may vary. Our free reference uses finite free convolution \cite{MSS,AFPU}. For two-resolvent eigenvector overlaps see \cite[eqs.~(6)--(11)]{BBP}; Theorem~\ref{thm:replica} uses independent samples for recovery. Gaussian skew-product densities and the isotropic Brownian area limit appear in \cite{FILZ,DO}.
\subsection{Proofs of the finite-dimensional identities}
\begin{proposition}[The two-spectrum mean for independent rows]\label{reb:mean-rows}
For a complex $d\times d$ matrix $Z$ the mean determinant is
\[
 \mathcal T(Z):=\E\det(Z-iA_d^X)
   =\left.\prod_j(1-\lambda_j^2\partial_y^2)\det(Z+yQ_d)\right|_{y=0}.
\]
It is nonzero for $\operatorname{Im} Z\succ0$. Its scalar restriction gives \eqref{eq:two-spectrum-polynomials}, and its derivative at $zI$ gives the root weights above.
\end{proposition}
\begin{proof}
First truncate the time space to finitely many skew pairs. For row and column sets of size $m$, independence, centering and isotropy give
$\E[\det X_{I,S}\det X_{J,T}]=m!\mathbf1_{I=J,S=T}$.
Cauchy--Binet applied twice yields
$\E[\wedge^m A_d^X]=m!\tr(\wedge^m K_d)\,\wedge^m Q_d$.
The time trace vanishes in odd degree and equals $e_k(\lambda^2)$ at $m=2k$. Expansion in minors of $Z$ gives the formula. Under spectral truncation the principal Pfaffians converge in $L^2$ by the minor-summation formula~\cite{IW} and the same orthogonality, so $\det(I+A_d^X)$ converges in $L^1$. Paired singular values give $|\det(Z-iA_d^X)|\le (1+\|Z\|_{\op}^2)^{d/2}\det(I+A_d^X)$, passing the identity to Hilbert--Schmidt kernels.

For $\operatorname{Im}Z\succ0$ all $y$-zeros of $\det(Z+yQ_d)$ lie in the open lower half-plane. Each $1\pm\lambda_j\partial_y$ preserves this property since the logarithmic derivative of a nonconstant such polynomial has negative imaginary part on the closed upper half-plane. Coefficient convergence and Hurwitz give the infinite product case, using $\mathcal T(iI)=i^dP_d(1)\ne0$. Thus $p_d(z)=\mathcal T(zI)$ is real-rooted. In any real orthonormal basis $G(x)=\E\det(I+tA_d^X\operatorname{diag}x)$ is nonzero for $\operatorname{Re}x_i>0$ by taking $Z=(i/t)(\operatorname{diag}x)^{-1}$; the principal minors have even degree.

Let $R(z)=\E\operatorname{adj}(zI-iA_d^X)/p_d(z)$. For each unit $u\in\mathbb C^d$, the rank-one identity $\mathcal T(zI+h uu^*)=p_d(z)(1+h\,u^*R(z)u)$ implies $\operatorname{Im}(u^*R(z)u)\le0$ for $\operatorname{Im}z>0$. This rational function is real off its real poles on the real axis and asymptotic to $1/z$. Its poles are simple with positive residues. Polarization yields the asserted matrix weights. Their traces follow from $\tr R=p_d'/p_d$. At $z=i/t$ the self-adjoint part of $zR(z)$ equals both $\int(1-w_t)\,\dd\Xi_d$ and $\E[D(I-W_d(t))]/\E D$. In a $Q_d$-eigenbasis the shift formula is invariant under sign changes and under rotations in repeated eigenspaces. Hence $R(z)$ is diagonal and scalar on each such eigenspace. Deleting coordinate $i$ gives $zR_{ii}(z)=p_d(z)|_{q_i=0}/p_d(z)$, since setting $q_i=0$ leaves a scalar $z$ block. At $z=i/t$ this is $P_d(t)|_{q_i=0}/P_d(t)=1-q_i\partial_{q_i}\log P_d(t)$ by multiaffinity.
\end{proof}

\begin{proof}[Proof of Theorem~\ref{reb:mode}]
As $v$ increases from $0$ to infinity, $vw=D_Q(v)$ increases from $0$ to $\rank Q_d$ and $w$ decreases from $\tr Q_d$ to $0$. Hence $D_Q(v)=D_K(tw)$ has a unique crossing. Inverting both occupations and using $vw=x$ gives \eqref{reb:inverse-scales}.

We need a ratio bound for a sum $S$ of independent Bernoulli trials with probabilities $0<p_i<1$, mean $a>0$ and masses $s_n$. On consecutive positive masses,
\[
 R_n:=\frac{n s_n}{a s_{n-1}}
 \quad\text{satisfies}\quad R_n\ge1\ (n\le a),\qquad R_n\le1\ (n\ge a+1).
\]
For the first bound, multiply the success odds by $z=n/a\le1$. The new mean is at least $za=n$. By Darroch's mode theorem \cite{Darroch} every mode lies within distance strictly less than one of the mean. Log-concavity of Bernoulli sums by Newton's inequalities then gives $z s_n/s_{n-1}\ge1$.

The upper bound is the Poisson--binomial mass-ratio inequality of D\"umbgen and Wellner \cite[proof of Proposition~1, (7)]{DW}.
Let $S,T$ be Bernoulli sums with positive success odds, masses $s_n,h_n$ and means $a,\ell a/2$, with $\ell=1$ or $2$. Suppose a random index $J$ has masses proportional to $(2k)!s_{2k}h_{\ell k}/a^{2k}$, and its generating function is real-rooted. Its consecutive mass ratio $R_{2k}R_{2k-1}h_{\ell k}/h_{\ell(k-1)}$ is, with $m=a/2$, at least one for $k\le\lfloor m\rfloor$ and at most one for $k\ge\lceil m\rceil+1$: Darroch and log-concavity place each of the $\ell$ factors in the last ratio on the correct side of one. The support contains $\lfloor m\rfloor$, so there is a mode in $[\lfloor m\rfloor,\lceil m\rceil]$. The real-rooted generating function factors into Bernoulli generating functions, and Darroch implies $|2\E J-a|<4$.

For a finite temporal list take success odds $vq_i,t^2w^2\lambda_r^2$, omit zeros, and use $a=\theta_d(t)$ and $\ell=1$. The means have the stipulated values, and Proposition~\ref{reb:mean-rows} identifies the generating function with $P_d(t\sqrt{s})/P_d(t)$, giving $2\E J=d\Theta_d(t)$. For infinite temporal rank the coefficients converge under truncation and $D_K(tw)$ converges uniformly for $0\le w\le\tr Q_d$ by square summability of $\lambda$. The unique crossings converge, proving \eqref{reb:two-scalar-accuracy}.

Both scalar error signs occur for $d=2$, $t=1$, spatial weights $R,R^{-1}$ and one time pair $\lambda=R^{1/2}$ or $R^{-1/2}$. Then $d\Theta_d=4\lambda^2/(1+2\lambda^2)$ tends to two or zero. At $v=\lambda$, $D_Q(v)\to1$ and $D_K(w)=2D_Q(v)^2/(1+D_Q(v)^2)\to1$, and the crossing lies between these two readings. Padding with zeros proves sharpness of the normalized order $1/d$.

In any real orthonormal basis the nonnegative principal minors give a subset law $\mu_1$ with generating function $G(z)/G(\mathbf1)$. We use the same tangent inequality here and in Proposition~\ref{short:necessity}. Let $F$ be homogeneous of degree $n>0$ with nonnegative coefficients, $F(\mathbf1)=1$, and zero-free when $|\arg z_j|<\pi/4$ for all $j$. By \cite[Corollary~2]{Anari26}, $F(\sqrt{x})^{2/n}$ is concave. Its tangent at $\mathbf1$ gives
\begin{equation}\label{eq:sector-tangent}
 F(z)^{2/n}\le n^{-1}\sum_j\partial_jF(\mathbf1)z_j^2 \qquad(z\ge0).
\end{equation}
For the subset indicators $B_i$, take $F(z,w)=\E_{\mu_1}\prod_i z_i^{B_i}w_i^{1-B_i}$. It has degree $d$ and equals $(\prod_i w_i)G((z_i/w_i)_i)/G(\mathbf1)$. All ratios have positive real part in the required sector. Taking $(z_i,w_i)=(f_i(1),f_i(0))$ with $\E_{\mu_1}f_i(B_i)^2=1$ makes the right side one. Rescaling gives $\|\prod_i f_i(B_i)\|_1\le\prod_i\|f_i(B_i)\|_2$. Here and below, for a joint law $\tau$ this product bound yields $\sum_i\mathrm{KL}(\eta_i\Vert\tau_i)\le2\mathrm{KL}(\eta\Vert\tau)$ by entropy duality with $f_i=e^{h_i/2}$ \cite{AJKPV}.

Given $A=A_d^X$, sample $S$ with probabilities $t^{|S|}\det(A_{S,S})/D$. Jacobi's formula gives inclusion probabilities $(I-(I+tA)^{-1})_{ii}=W_{ii}$, since $W$ is the symmetric part of $I-(I+tA)^{-1}$. The original law and its tilt by $D/P_d(t)$ thus induce subset laws $\eta$ and $\mu_1$ with marginals $r_i$ and $m_i$. Data processing gives $\mathrm{KL}(\eta\Vert\mu_1)\le\log P_d(t)-\E\log D=dJ_d(t)$. In an eigenbasis of $\E W-N$, Pinsker's inequality $2(r_i-m_i)^2\le\operatorname{kl}(r_i\Vert m_i)$ proves the first bound in \eqref{eq:finite-spatial-gap}.
At every positive field $x$ the indicator covariance $C$ is diagonally dominant. Indeed,
fix $i$ and set $z_j=w^{\epsilon_j}$, $\epsilon_j=\pm1$, for $j\ne i$ in
$G(x_1z_1,\ldots,x_dz_d)=a(w)+z_i b(w)$. By half-plane nonvanishing
$h=(a-b)/(a+b)$ maps the right half-plane to the closed unit disk. At
$1$, $h=1-2p_i$ and $h'=-2\sum_{j\ne i}\epsilon_j C_{ij}$, with
$p_i$ the inclusion probability. Schwarz--Pick gives
$|h'|\le(1-h^2)/2=2C_{ii}$. Taking signs gives $\sum_{j\ne i}|C_{ij}|\le C_{ii}$, hence $C\preceq2\operatorname{diag}(C_{ii})\preceq I/2$.

In a $Q_d$-eigenbasis the weights of $\mu_1$ are proportional to $b_{|S|}\prod_{i\in S}q_i$, with $b_{2k}=(2k)!t^{2k}e_k(\lambda^2)$ and $b_{2k+1}=0$. Thus the time kernel only reweights the cardinality. The marginals $m_i$ are the eigenvalues of $N_d$ by spatial differentiation. If the law charges nonempty sets, the empty set and every pair of positive weights have positive mass, so $0<m_i<1$ on the active coordinates. Swapping $i,j$ gives $q_jm_i(1-m_j)-q_i m_j(1-m_i)=(q_j-q_i)C_{ij}$, where $C_{ij}$ is their inclusion covariance. Thus for $r_j=m_j/[q_j(1-m_j)]$,
\[
 m_i-p(r_jq_i)=\frac{(q_j-q_i)C_{ij}}{q_j(1-m_j)+q_i m_j},\qquad
 \sum_i|m_i-p(r_jq_i)|\le1
\]
since the absolute coefficient of $C_{ij}$ is at most $1/\min(m_j,1-m_j)$ and $\sum_{i\ne j}|C_{ij}|\le m_j(1-m_j)$. The profiles at $\min_jr_j,\max_jr_j$ bracket $m$ coordinatewise and differ in $\ell^1$ by at most two. Therefore $\|N_d-p(vQ_d)\|_{\mathcal S_1}\le\max\{2,|\tr N_d-D_Q(v)|\}$ for any $v>0$: between the profiles the error is at most two, and outside them all differences have the same sign. This also covers $m=0$. Taking $v=v_d(t)$ and applying \eqref{reb:two-scalar-accuracy} proves the matrix bound.
For rank-one $Q_d$ the sample is zero, whereas $\theta_d(t)\uparrow1$ for any nonzero $K_d$. Thus a uniform $o(1/d)$ remainder is impossible even at zero Jensen gap.
\end{proof}

\subsection{The free reference and its spatial resolvent}\label{sec:free-reference}
Write $a_j=d\lambda_j$, $\rho_d=d^{-1}\sum_i\delta_{q_i}$, and let
$\pi_d$ be the symmetric freely infinitely divisible law with Voiculescu transform
\begin{equation}\label{eq:free-time-law}
 \varphi_d(z)=\frac2d\sum_j\frac{z a_j^2}{z^2-a_j^2}.
\end{equation}
This free L\'evy--Khintchine law \cite{BP,BG} gives jumps $\pm a_j$ of rate $1/d$ and folded canonical measure $\eta_d=(2/d)\sum_j a_j^2(1+a_j^2)^{-1}\delta_{a_j}$. Square summability justifies the series. Its even cumulants $(2/d)\sum_j a_j^{2k}$ grow at most exponentially and its odd cumulants vanish, hence the moment--cumulant formula gives bounded support. Let $Y\sim\pi_d$ be free from $\mathcal A_d=M_d(\mathbb C)$ with $\tr_d=d^{-1}\tr$, $\tau$ the joint trace and $Z=Q_d^{1/2}YQ_d^{1/2}$; write $\mathbb E_{\mathcal A_d}$ for conditional expectation.
\begin{proposition}[Finite and free representatives]\label{prop:free-reference}
With $\chi_Q(x)=\det(xI-Q_d)$ and $\partial=\partial_x$, the mean characteristic polynomial is
\begin{equation}\label{eq:finite-free-factor}
 p_d=\chi_Q\boxtimes_d\left[\prod_j(1-\lambda_j^2\partial^2)x^d\right].
\end{equation}
The polynomial in brackets is the finite free additive convolution of the individual polynomials $(1-\lambda_j^2\partial^2)x^d$. The spectral law of $Z$ on $\mathbb R$ is the ordinary free multiplicative convolution $\rho_d\boxtimes\pi_d$. For every $t>0$, with $v$ from \eqref{reb:two-spectral-equation},
\begin{equation}\label{eq:free-marked-resolvent}
 \mathbb E_{\mathcal A_d}\big[(I+itZ)^{-1}\big]=(I+vQ_d)^{-1},\qquad
 \tau\bigl(w_t(|Z|)\bigr)=\theta_d(t)/d.
\end{equation}
The standard squaring identity for a symmetric free multiplicative factor~\cite{AP09} gives $\mathcal L(Z^2)=\rho_d\boxtimes\rho_d\boxtimes\mathcal L(Y^2)$. Taking real parts yields $\mathbb E_{\mathcal A_d}[w_t(|Z|)]=M_d(t)$, the explicit approximation in Theorem~\ref{reb:mode}. If an input vanishes, use $v=0$ in \eqref{eq:free-marked-resolvent}.
\end{proposition}
\begin{proof}
Write $(d)_n=d!/(d-n)!$. For $f(x)=\sum_n(-1)^ne_n(f)x^{d-n}$, finite free multiplication divides $e_n(f)e_n(g)$ by $\binom dn$, while finite free addition multiplies the differential symbols acting on $x^d$ \cite{MSS}. Thus the bracket coefficient at $2k$ equals $(-1)^k(d)_{2k}e_k(\lambda^2)$; $(d)_{2k}/\binom d{2k}=(2k)!$ gives \eqref{eq:two-spectrum-polynomials}, with infinite products taken coefficientwise.

Let $M(z)=\mathbb E_{\mathcal A_d}[(zI-Z)^{-1}]$, $g=\tr_d(Q_dM)$ and $R_Y(w)=\varphi_d(1/w)$. By scalar freeness and \cite[Theorem~3.5, eq.~(3.6)]{NSS02}, inserting the spatial factors gives
$\kappa_n^{\mathcal A_d}(ZB_1,\ldots,ZB_{n-1},Z)=Q_d\kappa_n(Y)\prod_{j<n}\tr_d(Q_dB_j)$.
The matrix $R$-transform is therefore $B\mapsto Q_dR_Y(\tr_d(Q_dB))$, and the moment--cumulant recursion gives
\begin{equation}\label{eq:free-dyson}
 [zI-Q_dR_Y(g(z))]M(z)=I .
\end{equation}
Initially valid near infinity by boundedness, this holds on $\mathbb C_+$ by analyticity, since $R_Y(w)=(2/d)\sum_j a_j^2w/(1-a_j^2w^2)$ converges locally uniformly on $\mathbb C_-$ and $-\operatorname{Im}g=(\operatorname{Im}z)\tau(Q_d^{1/2}|zI-Z|^{-2}Q_d^{1/2})>0$ for $Q_d\ne0$.

At $z=i/t$, symmetry gives $g=-iy$ and $R_Y(-iy)=-iG_d(y)/y$, $G_d(y)=(2/d)\sum_jw_y(a_j)$. Hence $M=-it(I+vQ_d)^{-1}$, $v=tG_d(y)/y$. Taking the $Q_d$-weighted trace yields $y=t\tr_d[Q_d(I+vQ_d)^{-1}]$. Hence $vy/t=G_d(y)=\theta_d(t)/d$ by uniqueness, and multiplication by $i/t$ gives \eqref{eq:free-marked-resolvent}.
\end{proof}

\section{The row criterion}\label{reb:universality}
Write $B=Q^{1/2}$ and $F(Z)=d^{-1}\log\det(I+tBZKZ^*B)$.

\begin{proposition}[Comparison on bounded windows]\label{thm:q-transfer}
Assume \eqref{eq:pc}. If $\|Q\|_{\op}\le1$, $d\|K\|_{\HS}^2\le E$ and $d\|K\|_{\op}\le M$, then at each fixed $t>0$
\[
 \E F(X)-d^{-1}\log P_d(t)\longrightarrow0,\qquad
 d^{-1}\|\E W_d(t)-M_d(t)\|_{\mathcal S_1}\longrightarrow0 .
\]
Rescaling $K$ permits any fixed bound on $\|Q\|_{\op}$. Taking traces gives comparison with $\Theta_d$ by \eqref{reb:two-scalar-accuracy}.
\end{proposition}
\begin{proof}
In a finite skew basis take reference rows $g_i$ with independent uniform coordinates on $[-\sqrt3,\sqrt3]$. Their product law $m$ satisfies $W_2^2(\nu,m)\le 2C\mathrm{KL}(\nu\Vert m)$ with $C=3$: a conditional quantile differs from the uniform quantile by at most $2\sqrt3$ times total variation, so apply Pinsker and the entropy chain rule. Write $f=dF$, $b=\|Q\|_{\op}$, $h=\|K\|_{\HS}$, $k=\|K\|_{\op}$. Since $R=(I+tBZKZ^*B)^{-1}$ satisfies $\|R-R^*\|_{\op}\le1$, differentiation gives $\|\nabla f(Z)\|_F\le tb\|ZK\|_F$. Couple $Z\sim m$ optimally to $Y\sim e^fm/P_d$ and put $H=\mathrm{KL}(\mathcal L(Y)\Vert m)$. Integration along their segment and $\E\|ZK\|_F^2=dh^2$ give $dJ_d\le tbh\sqrt{2CdH}-(1-Ctbk)H$. Thus $J_d\le Ct^2b^2h^2/[2(1-Ctbk)]$ for $Ctbk<1$, hence $O(d^{-1})$ in the window. Gaussian products have $C=1$ \cite{GL}.

Let $k_r(S)$ sum the top $r=\lceil ad\rceil$ eigenvalues of $S\succeq0$ (pad by zeros). The supremum $D_d(a)$ in \eqref{eq:pc-stat} over all finite ranks $p\ge ad$ tends to zero, since ranks above $d/4$ split into ranks in $[d/10,d/4]$ for large $d$. Integrating over the spectral projections of $S$, with the bound twice the rank below $r$, gives
\begin{equation}\label{short:layers}
 d^{-1}\sum_i\E|X_i^*SX_i-\tr S|\le2k_r(S)+D_d(a)\tr S .
\end{equation}
For an average over $n$ labels the high-rank coefficient becomes $\min(2,dD_d(a)/n)$. For the reference rows $D_d(a)\le\sqrt{2/r}$ by fourth moments.

Take the families independent and put $L=\sum_i K^*(X_iX_i^*+g_ig_i^*)K$, so $\E\tr L\le2E$. For either family $z_i$ set $h_i=d\|K^*z_i\|^2$. Applying \eqref{short:layers} to $dKK^*$, whose norm is at most $M^2/d$, gives $d^{-1}\sum_i\E|h_i-dh^2|\to0$ by taking $a\downarrow0$ after $d\to\infty$. Discarding $h_i>2E+1$ costs $o(1)$ in expected trace. Each retained sum $U$ satisfies $\E\tr U^2\le(2E+1)E/d+d(d-1)\tr((K^*K)^2)=O(d^{-1})$ by truncation on the diagonal and independence off it. Hence $k_r(U)\le\sqrt{r\tr U^2}$ and subadditivity give $\lim_{a\downarrow0}\limsup_d\E k_r(L)=0$.

Independently take $Z_i=X_i$ with probability $s$ and $Z_i=g_i$ otherwise. All these mixtures have the same mean determinant. Replacing row $i$ by $x$ gives $Z^{i:x}$. Differentiating the product mixture yields $\phi'(s)=\sum_i\E[F(Z^{i:X_i'})-F(Z^{i:g_i'})]$ for $\phi(s)=\E F(Z)$, with fresh candidates independent of $Z$. Put $Y=BZ$, $S=Y^*Y$, $A=YKY^*$. The symmetric part of $(I+tKS)^{-1}tK$ is $T(S)=t^2K^*Y^*(I+t^2A^*A)^{-1}YK\preceq t^2 L$. For $q=Q_{ii}>0$ set $S_-=Z^*(Q-Qe_ie_i^*Q/q)Z$ and $\ell=q^{-1}\sum_{j\ne i}Q_{ij}Z_j$. Both omit row $i$. Completing the square adds $q(x+\ell)(x+\ell)^*$ to $S_-$, so the determinant lemma gives the row factor $1+\|c+Vx\|^2$, with $V=(qT(S_-))^{1/2}$ and $c=V\ell$. Rows with $q=0$ give zero increment. The unsymmetrized resolvent difference has rank at most one. Since $T(S_-),T(S)\preceq t^2L$ and $q\le1$, we have $\|V^*V-qT(S)\|_{\mathcal S_1}\le2t^2\|L\|_{\op}$.

For a centered isotropic candidate set $C=V^*V$, $z=\langle c/\|c\|,Vx\rangle$ (zero if $c=0$). In $\log(1+x^*Cx-z^2+(\|c\|+z)^2)$ replacing the nonnegative $x^*Cx-z^2$ by $\tr C$ costs at most $\E|x^*Cx-\tr C|+\E z^2$ in expectation. The resulting function of $z$ has second derivative bounded by two. Since $\E z=0$,
\[
 0\le\log(1+\|c\|^2+\tr C)-\E\log(1+\|c+Vx\|^2)
 \le\E|x^*Cx-\tr C|+2\|V\|_{\op}^2.
\]
Condition on the background. Both candidates have the same term $\log(1+\|c\|^2+\tr C)$. In the deviation, replacing $C$ by $qT(S)$ costs at most $2\|C-qT(S)\|_{\mathcal S_1}$ by isotropy. Now $T(S)$ is common to all labels. Applying \eqref{short:layers} to both families yields, uniformly in $s$,
$|\phi'(s)|\le t^2\E[16k_r(L)+2\max(D_d(a),\sqrt{2/r})\tr L]$, since $\|L\|_{\op}\le k_r(L)$. Integrating gives $\E F(X)-\E F(g)\to0$. Time truncation preserves the projection bounds and converges entrywise in $L^2$. The $t/\sqrt d$ Frobenius Lipschitz bound for $F$ as a function of the skew matrix and \eqref{eq:two-spectrum-polynomials} pass the estimates to Hilbert--Schmidt $K$. Apply \eqref{eq:finite-spatial-gap} for the matrix assertion.
\end{proof}

The logarithmic gap can be recovered from the ordinary response. Let $\delta=\sup_{0<s\le t}|\Theta_d(s)-\E\tr_d W_d(s)|$ and let $m_2=[(\tr Q)^2-\tr(Q^2)]\|K\|_{\HS}^2/d$ be the common second moment of the root and mean sample laws, by the degree-two coefficient. Since both responses lie in $[0,m_2s^2]$, logarithmic differentiation gives, for $\delta>0$,
\[
 J_d(t)=\int_0^t(\Theta_d(s)-\E\tr_d W_d(s))\,\frac{\dd s}{s}
 \le \frac{\delta}{2}\left[1+\log_+(t^2m_2/\delta)\right].
\]
The bound is zero for $\delta=0$. Thus tracking implies $J_d(t_d)\to0$ for bounded $t_d^2m_2$. For sharpness take independent sign coordinates, $Q$ the first two-coordinate projection and $K$ one unit skew pair. Then $D=1,1+4t^2$ equiprobably, $m_2=4/d$ and the response error increases to $1/d$. In dimension two the active marginals are $r=2t^2/(1+4t^2)$ and $m=2t^2/(1+2t^2)$. As $t\to\infty$, $\operatorname{kl}(r\Vert m)/(dJ_d(t))\to1$, proving sharpness of the entropy factor two in Theorem~\ref{reb:mode}. At $t=e^d$ the gap and its integral upper bound both tend to one; these rows satisfy \eqref{eq:pc} since $\Var(\|PX_i\|^2)\le2\rank P$.

\begin{proposition}[Determinant tilt and the one-rank row test]\label{short:necessity}
Let $Q=I$, $t>0$ and $K$ be skew-adjoint and Hilbert--Schmidt. Write $\mu=\bigotimes_i\mu_i$ for the row law, $\tau=D\mu/P_d(t)$ for its tilt and $\tau_i$ for the tilted marginals. Let $E_j$ project onto the skew plane with amplitude $\lambda_j$, and set $b_j=(t\lambda_j)^2$. From \eqref{eq:two-spectrum-polynomials} define
\[
 \gamma_j=b_j\partial_{b_j}\log P_d(t),\qquad B_t=\frac1d\sum_j\gamma_j E_j,\qquad
 \theta=\tr B_t=\Theta_d(t).
\]
Then $\dd\tau_i=(1-\theta+X_i^*B_tX_i)\,\dd\mu_i$, and $\sum_i\mathrm{KL}(\eta_i\Vert\tau_i)\le2\,\mathrm{KL}(\eta\Vert\tau)$ for every joint row law $\eta$ with marginals $\eta_i$.
For $\theta>0$ the normalized temporal energy $R_i^*=X_i^*B_tX_i/\theta$ satisfies
\[
 J_d(t)\ge\frac1{2d}\sum_i\mathrm{KL}(\mu_i\Vert\tau_i)
 \ge\frac{\theta^2}{4d}\sum_i\E(\sqrt{R_i^*}-1)^2 .
\]
When $\theta=0$ the tilt is trivial.

For $K=J_P/d$, $t=1$, let $P$ have positive even rank $p\le d/4$. Put $x_i=PX_i$, $R_i=\|x_i\|^2/p$, $H_i=\E(\sqrt{R_i}-1)^2$ and $\rho=p/d$. Write $H=d^{-1}\sum_iH_i$, $e_i=\rho\E|R_i-1|$ and $\varepsilon_P=d^{-1}\sum_i e_i$.
For every $J_P$, $R_i^*=R_i$ and $J_d(1)\ge\theta^2H/4\ge\rho^2H/31$.
Let $J^\circ$ be the gap for independent rows of length $\sqrt p$ with independent fair signs and directions sampled with weight $R_i$ from the original laws. Then $|J_d(1)-J^\circ|\le(1-1/d)\varepsilon_P/2\le\rho\sqrt H$.
Thus $J_d(1)\to0$ on a sequence of planes exactly when $\rho^2H\to0$ and $J^\circ\to0$. If $x_i=\sqrt{R_i}\,u_i$ with $R_i$ independent of $u_i$, $\|u_i\|=\sqrt p$, then $0\le J_d(1)-J^\circ\le2(1-1/d)\rho H$. At rank zero take $H=e_i=\varepsilon_P=J^\circ=0$.
For the prescribed rank $r_d$, \eqref{eq:pc} is equivalent to $\sup_{\rank P=r_d}J_{d,P}\to0$, or to $F(X)-d^{-1}\log P_d(1)\to0$ in probability for every sequence of rank-$r_d$ planes with the chosen complex structures. The function $h_d(u)=d^{-1}\sum_i\E|X_i(u)|$ alone determines whether these equivalent conditions hold.
\end{proposition}
\begin{proof}
First let $K$ have finite rank and project the rows to its range, writing $x_i$ for their projections. For bounded nonnegative $g_i$ on the full rows, set $a_i=\E g_i$, $C_i=\E[g_ix_ix_i^*]$. Each squared principal Pfaffian is quadratic in its selected rows, and independent rows enter each principal minor symmetrically, so $\int\prod_i g_i\,\dd\tau$ is a symmetric $d$-linear form $\mathcal Z$ in $(a_i,C_i)$. Write $G(a,C)$ for its diagonal. In a $C$-eigenbasis with eigenvalues $c_h$ and $T=tK$,
\[
 P_d(t)G(a,C)=\sum_{|I|\le d}(d)_{|I|}\det T_{I,I}\,a^{d-|I|}\prod_{h\in I}c_h .
\]
This follows from Proposition~\ref{reb:mean-rows} applied to $C^{1/2}TC^{1/2}$ for $C\succeq0$ and polynomial continuation. All coefficients are nonnegative.

For $M\ge d$ the size-$d$ principal-minor polynomial of $T\oplus I_M/M$ has the sector in \eqref{eq:sector-tangent} by \cite[Lemma~12]{AASV}. Identify its auxiliary variables with $a$. Multiplying by $d!/P_d(t)$ gives coefficientwise convergence to $G$, since $\binom{M}{m}/M^m\to1/m!$ for $m=d-|I|$. After normalization at one, \eqref{eq:sector-tangent} passes to $G$ in every basis.

At $e=(1,I)$ we have $G(e)=1$. Separate rotations of the skew pairs make the derivative in $C$ block diagonal and scalar on each plane, and varying one coordinate of pair $j$ gives $\gamma_j$. Hence $G'_C(e)[H]=d\tr(B_tH)$ and homogeneity gives $G'_a(e)=d(1-\theta)$. Positivity gives $|G(a,C)|\le G(|a|,|C|)$. The tangent bound at the all-ones vector now reads
\[
 |G(a,C)|^{2/d}\le(1-\theta)a^2+\tr(B_tC^2)=:\|(a,C)\|^2.
\]
The constant and quadratic terms give $\theta<1$ and $B_t>0$ on the range when $d\ge2$, $K\ne0$. Banach's diagonal norm theorem for real symmetric multilinear forms \cite{Banach} gives $|\mathcal Z|\le\prod_i\|(a_i,C_i)\|$. With all but one slot at $e$, the form equals $G'(e)/d$, giving the marginal density $u_i=1-\theta+X_i^*B_tX_i$. Compressing multiplication by $g_i$ via the isometry $v\mapsto\langle x_i,v\rangle$ gives $C_i^2\preceq\E[g_i^2x_ix_i^*]$. Together with $a_i^2\le\E g_i^2$ this yields
\[
 \int\prod_i g_i\,\dd\tau\le\prod_i\left(\int g_i^2\,\dd\tau_i\right)^{1/2}.
\]
For $K=0$ or $d=1$, $\tau=\mu$. Under spectral truncation $D$ converges in $L^1$ by Proposition~\ref{reb:mean-rows}, so all integrals in the product bound converge. The $\gamma_j$ and their sum converge by \eqref{eq:two-spectrum-polynomials}, giving trace-norm convergence of $B_t$ and hence $L^1$ convergence of its row energies by isotropy. Entropy duality now gives the assertion for $\eta$. For $\eta=\mu$, the entropies are $dJ_d(t)$ and $\E[u_i-1-\log u_i]$. Since $u-1-\log u\ge(u-1)^2/[2\max(1,u)]$ and $\max(1,u_i)\le(1+\sqrt{R_i^*})^2$ when $\theta>0$, the energy bound follows.

On a flat plane symmetry gives $B_t=\theta P/p$. In $P_d(\sqrt z)/P_d(1)$ the coefficient law $w_k\propto (d)_{2k}\binom{p/2}k/d^{2k}$ has mean $d\theta/2$. Since $(k+1)w_{k+1}\ge(9/16)(p/2-k)w_k$, the mean is at least $9p/50$. For $p=0$, $D=1$. The same multilinear bound yields the fixed-norm maximum of each squared principal Pfaffian in expectation. Take $d=2k\le p$, $a_i=0$ and let $t\to\infty$. Then $\theta\to1$ and equality holds at $C_i=\|C_i\|_{\HS}I/\sqrt p$. For $p>0$ put $\nu_i=R_i\mu_i$.
Write $x_i=\sqrt{pR_i}\,\theta_i$ with arbitrary unit direction at zero. The fixed-length replacement is isotropic since $p\int\theta_i\theta_i^*\,\dd\nu_i=I$. Replacing row $i$ changes the mean log factor by $\E g(R_i)$ up to sign, where $y=p\theta_i^*M_i\theta_i$ and $g(r)=\log(1+ry)-r\log(1+y)$. Concavity and $|g'|\le y$ place $g(r)$ between zero and $y(1-r)$ in either order. Independence gives $\E M_i\preceq(d-1)I/d^2$. Since $\E(1-R_i)_+=\E(R_i-1)_+$, each sign contributes at most $(1-1/d)e_i/2$. Cauchy--Schwarz gives $e_i\le2\rho\sqrt{H_i}$. For $s=\sqrt{R_i}$ independent of $y$, Jensen gives a nonnegative gap increase. As $k(s)=\log(1+ys^2)$ has $|k''|\le2y$ and $\E(1-s)=H_i/2$, Taylor's bound gives $k(1)-\E[k(s)\mid y]\le[y/(1+y)+y]H_i$. Use $\E y\le(1-1/d)\rho$.

For $0<p\le r=r_d$ choose rank-$r$ $E\perp P$ and a uniform coordinate $(r-p)$-plane $E'\subset E$. The identity $P=\E(P+E')-(1-p/r)E$ gives $\varepsilon_P\le2\sup_{\rank L=r}\varepsilon_L$ by subadditivity of the centered quadratic-form deviation. Ranks up to $d/4$ split into at most four such blocks when $r>0$. If $r=0$ and $d\ge8$, the time space has dimension at most three. Since $\varepsilon_P\le2p/d$, this proves \eqref{eq:pc} equivalent to $\sup_{\rank P=r_d}\varepsilon_P\to0$. The lower bound and Proposition~\ref{thm:q-transfer} give the gap test.

To read the row condition from first moments, let $\sigma_P$ be uniform probability on the unit sphere of $P$ and $a_p=\int|u_1|\,\dd\sigma_{\mathbb R^p}(u)$. Spherical averaging of $|\langle x_i,u\rangle|$ and isotropy give
\[
 \frac1{a_p\sqrt p}\int h_d(u)\,\dd\sigma_P(u)
   =\frac1d\sum_i\frac{\E\|x_i\|}{\sqrt p}=1-H/2 .
\]
Now $\rho H\le\varepsilon_P\le2\rho\sqrt H$ gives the equivalent test $\sup_{\rank P=r_d}\rho H\to0$ (zero at rank zero). Equality of $h_d$ is zonoid equivalence for the uniformly chosen rows on finite subspaces \cite{MSSZ}; it therefore preserves \eqref{eq:pc}.

For $Q=I$, deleting row $i$ changes $F$ by $d^{-1}\log(1+X_i^*M_iX_i)$ with $\E\tr M_i\le t^2(d-1)\|K\|_{\HS}^2$. Variance tensorization and concavity of $\log^2(e+u)$ give $\Var F\le d^{-1}\log^2(e+t^2(d-1)\|K\|_{\HS}^2)$. This tends to zero on the flat tests, proving the probability assertion. There $m_2=(d-1)p/d^2\le1/4$. Tracking the mean ordinary response at rational $0<t<1$ extends to this interval by $|\partial_{\log t} w_t|\le1/2$. The integrand for $J_d(1)$ is bounded in absolute value by $m_2t$. Dominated convergence gives necessity from these flat ordinary responses alone, whose reference depends only on $d,r_d$.
\end{proof}
For $p=d/4$ even, take $x_i=\sqrt{R_i}\,u_i$ with independent sign coordinates in $u_i$, independent of the radii, so \eqref{eq:direct-jensen} gives $J^\circ=O(d^{-1})$. For fair $R_i=1\pm a_i$, $0\le a_i\le1/2$, we have $H_i\asymp a_i^2$. Either $a_i=d^{-1/4}$ throughout or $a_i=1/2$ on $\lfloor\sqrt d\rfloor$ rows and zero elsewhere gives a gap of order $d^{-1/2}$. Killing each row with probability $\delta=d^{-1/2}$ and rescaling survivors by $(1-\delta)^{-1/2}$ gives $H_i=2(1-\sqrt{1-\delta})\asymp\delta$, $e_i^2\asymp\delta^2$ and a gap of order $\delta$.

For $x_i$ uniform on the signed coordinate axes of length $\sqrt p$, pair consecutive axes by $J$. A pair receiving $n$ rows, $b$ on its first axis, contributes $1+b(n-b)/16$ to the determinant product. Given the totals, the splits are independent fair binomials. Each pair with $n=2$ contributes $\log(33/32)-\tfrac12\log(17/16)>0$ to the conditional gap. Their expected count divided by $d$ tends to $4e^{-8}$ since $n\sim\operatorname{Bin}(d,2/p)$. Thus $\liminf J_d(1)>0$ even at $H=0$. Placing the two spherical laws independently on orthogonal planes explains the supremum over subspaces.
\section{Spatial comparison at all scales}\label{sec:unrestricted}

\subsection{Bounded spatial spectra}
For $H=iA_d^X$ the frame $\binom{I}{tH}$ has volume $\det(I+t^2H^2)^{1/2}=\det(I+tA_d^X)$ by paired eigenvalues. For every Hermitian $H$ its graph projection has diagonal blocks $I-w_t(|H|)$ and $w_t(|H|)$, by the graph frame $\binom{I}{tH}(I+t^2H^2)^{-1/2}$. If $\rank(H-H')\le k$, the two graphs intersect in dimension at least $d-k$. Their projections differ by a contraction of rank at most $2k$, so pinching to the two diagonal blocks gives $\|w_t(|H|)-w_t(|H'|)\|_{\mathcal S_1}\le k$.
Put $x_j=d\lambda_j$, $\rho_d=d^{-1}\sum_i\delta_{q_i}$ and
\[
 \eta_d=\frac2d\sum_j\frac{x_j^2}{1+x_j^2}\delta_{x_j},\qquad
 g_y(s)=\frac{y^2(1+s^2)}{1+y^2s^2},\qquad
 H(v)=\int\frac{vq}{1+vq}\,\rho(\dd q),\quad L(v)=H(v)/v .
\]
Write $G_\eta(y)=\int g_y\,\dd\eta$, $G_d=G_{\eta_d}$, and $H_d,L_d$ for the spatial transforms of $\rho_d$. Set $g_y(\infty)=1$ for $y>0$. The finite scalar crossing reads $\theta_d(t)/d=H_d(v)=G_d(tL_d(v))$.

\begin{lemma}[Mean matrix response below the spatial rank]\label{lem:marked-subcritical}\label{full:color}
Assume \eqref{eq:pc}, $\|Q_d\|_{\op}\le C$, $\rho_d\Rightarrow\rho$, and $\eta_d\Rightarrow\eta+\alpha\delta_\infty$ on $I_*$ with $\eta$ carried by $[0,\infty)$ and $\alpha<r=1-\rho(\{0\})$. Then for $t>0$,
\[
 d^{-1}\|\E w_t(|A_d^X|)-p(vQ_d)\|_{\mathcal S_1}\longrightarrow0,
\]
where $v=0$ for the zero source and otherwise $v>0$ is the unique solution of
\begin{equation}\label{eq:known-spatial-inverse}
 H(v)=\alpha+G_\eta(tL(v)).
\end{equation}
In particular the full response has limit $H(v)$.
\end{lemma}
\begin{proof}
If $d\lambda_j$ are bounded, then $d\|K_d\|_{\HS}^2=\int(1+s^2)\,\dd\eta_d$ is bounded and $\alpha=0$. Proposition~\ref{thm:q-transfer} gives normalized trace-norm convergence. Monotone squeezing in the scalar crossing gives $v_d(t)\to v$, including the zero source. The functions $p(v_dq)$ therefore converge uniformly on $[0,C]$ to $p(vq)$.
For unbounded amplitudes the number of temporal directions divided by $d$ on $(R,\infty]$ is $\int_{(R,\infty]}(1+s^{-2})\,\eta_d(\dd s)$. Choose continuity thresholds $a_d\uparrow\infty$ slowly enough that $n_d/d\to\alpha$, where $n_d$ counts directions above $a_d$. Delete $R<d\lambda_j\le a_d$: their limiting rank fraction is at most $(1+R^{-2})\eta((R,\infty))$ at any continuity point $R>1$. Let $Z$ collect the $n_d$ remaining high coordinates and put $W=Q_d^{1/2}Z/\sqrt d$. If $\alpha=0$ these high coordinates can also be deleted. For $0<\alpha<r$ we claim
\[
 \lim_{e\downarrow0}\limsup_d d^{-1}\E\#\{j\le n_d:s_j(W)\le e\}=0.
\]
Test on those coordinates with a skew isometry $J$. The bounded-window result for $WJW^*$ has full response $H(v)=\alpha(tL(v))^2/(1+(tL(v))^2)$. Since $H(v)\le\alpha<r$, $v$ stays bounded as $t\to\infty$ and $L(v)$ stays positive. Hence the limiting mean response tends to $\alpha$. Write $N_W(u)=\#\{j:s_j(W)>u\}$. By the singular-value product inequality, $\tr_d w_t(|WJW^*|)\le w_t(De)+[N_W(D)+N_W(e)]/d$. Use $\E N_W(D)\le Cn_d/D^2$ and take $d\to\infty$, $e\downarrow0$, $D\to\infty$, then $t\to\infty$.

Replace the high amplitudes by $M$, preserving their orientations. In both matrices raise the singular values of $W$ below $e$ to $e$ in a common thin singular-value decomposition, completing the left frame at zeros. This costs twice the exceptional count in rank in each matrix. The high Hermitian blocks now have common range $V$ and least singular values there at least $e^2 a_d$ and $e^2 M$, respectively. Their shared low block $B$ has $d^{-1}\E\tr B^2\le C^2(1+R^2)\sup_d\eta_d(I_*)$, so clipping it at $\pm D$ costs $O_{C,R}(D^{-2})$ in expected rank fraction. For either high block $A=A_V\oplus0$ with least singular value $b$ on $V$, compare the graph of $t(A+B)$ with
\[
 E_\infty=\{(y,z+tP_{V^\perp}By): y\in V^\perp,\ z\in V\}.
\]
For each displayed vector set $u=(tA_V)^{-1}(z-tP_VBy)$. The graph vector at $y+u$ differs from it by $(u,tBu)$, of norm at most $(1+t^2D^2)/(tb)$ times its norm. Thus the graph projection is within $(1+t^2D^2)/(tb)$ of $P_{E_\infty}$ in operator norm, since the planes have equal dimensions. Both high blocks have the same limit plane, so their vertical responses satisfy the corresponding comparison.
For fixed $R,M$ the capped array lies in the window already proved. Its temporal transform converges to $G_{\eta|_{[0,R]}}(y)+\alpha(My)^2/(1+(My)^2)$. The crossings converge to \eqref{eq:known-spatial-inverse} by monotone squeezing at fixed $v>0$, first as $M\to\infty$, then as $R\to\infty$. In the error bounds take successively $d\to\infty$, $M\to\infty$, $D\to\infty$, $e\downarrow0$ and $R\to\infty$ through continuity points. The triangle inequality and uniform continuity of $p(vq)$ for bounded $q$ conclude the proof.
\end{proof}
\subsection{Removing the spatial bounds}

\begin{lemma}[Spatial rescaling]\label{lem:uncap}
Let $H=H^*$, $D\succeq0$ of size $d$, $W=w_t(|H|)$ and $W'=w_t(|DHD|)$. Put
\[
 b=\tr_d\big[(1-D^{-2})_+(I-W)+(1-D^2)_+W\big],
\]
taking the first multiplier to be zero on $\ker D$. Then
\begin{equation}\label{eq:spatial-tracenorm}
 d^{-1}\|W'-W\|_{\mathcal S_1}\le2\sqrt b .
\end{equation}
For $D\succeq I$ the normalized trace gain lies in $[0,2b]$. For $D\preceq I$ the same holds for the loss. The scalar factor two is sharp.
\end{lemma}
\begin{proof}
Write $P$ for the graph projection. For invertible $D$ the new graph is $T\operatorname{ran}P$, $T=D^{-1}\oplus D$. Write $T=U^{-1}C$ with $C=T\wedge I$ and $U=T^{-1}\wedge I$. Their contracted coordinates are disjoint.

For any $d$-plane $E$ let $P_1$ project onto $CE$ and put $k=\tr[(I-C^2)P_E]$. Since $(I-P_1)CP_E=0$ and $(I-C)^2\preceq I-C^2$, we have $\|P_1-P_E\|_F^2=2\|(I-P_1)P_E\|_F^2\le2k$.
Also $CP_EC$ is a positive contraction supported on $CE$, so $P_1-CP_EC\succeq0$ has trace $k$. On $\ker(I-C)$ this means the projection increases in compression, with trace gain at most $k$.

Apply this first to $C$ on the original graph, then to $U$ on its intermediate orthogonal complement. The first cost is $db$. The second is at most $db$ because $I-U^2$ is supported on coordinates untouched by $C$, where the graph projection has increased. The resulting complement equals $(T\operatorname{ran}P)^\perp$. The two projection changes thus have Frobenius norm at most $\sqrt{2db}$ each. The diagonal blocks of their total change are $\pm(W'-W)$, giving $\|W'-W\|_F^2\le4db$ and hence \eqref{eq:spatial-tracenorm}. If $D\succeq I$, the first step increases the vertical trace of the graph, and the second increases the horizontal trace of its complement, each by at most $db$. These give the asserted vertical trace gain. Interchanging horizontal and vertical proves the contraction case. Continuity includes singular $D$. For $t=1$, $H=iu\left(\begin{smallmatrix}0&1\\-1&0\end{smallmatrix}\right)$ and $D=\operatorname{diag}(\sqrt s,1)$, $s>1$, the gain divided by $b$ is $2su^2/(1+su^2)\to2$.
\end{proof}
\begin{proof}[Proof of Theorem~\ref{thm:main}, sufficiency]
By \eqref{eq:finite-spatial-gap} we may use $M_d$ throughout.
It suffices to compare the mean matrices along arbitrary deterministic scales. Absorb the scale into $K_d$ and first take a subsequence with $x_d=\theta_d(1)/d\to x\in(0,1)$. At the scalar crossing $v_d$ replace $(Q_d,K_d)$ by $(v_dQ_d,K_d/v_d)$. This leaves the matrix unchanged and gives $H_d(1)=x_d=G_d(x_d)$. The time sources now have bounded mass since $g_y\ge\min(1,y^2)$. Extract spatial and temporal limits on $I_*$, denoted by $\sigma$ and $\eta+\alpha\delta_\infty$. Continuity gives
\[
 x=\int_{I_*}\frac{q}{1+q}\,\sigma(\dd q)=\alpha+G_\eta(x),
 \qquad \alpha\le x .
\]
Here $q/(1+q)=1$ at infinity. For $C>1$ clip $q$ by the continuous map $q_C=(q\vee C^{-1})\wedge C$. Write $H_C,L_C$ for the clipped transforms. The clipped limiting law has positive rank fraction one, and $\alpha\le x<1$, so Lemma~\ref{lem:marked-subcritical} applies. Its crossing $x_C=H_C(v_C)$ lies between the two sides evaluated at $v=1$, namely $H_C(1)$ and $\alpha+G_\eta(H_C(1))$. Both tend to $x$.

Write $W_C$ for the clipped response, $Q_{d,C}$ for the clipped covariance,
and $\|T\|_{1,d}=d^{-1}\|T\|_{\mathcal S_1}$.
Lemma~\ref{lem:marked-subcritical} gives $\|\E W_C-p(v_CQ_{d,C})\|_{1,d}\to0$.
Since $p(v_Cq_C)-p(q_C)$ has constant sign and
$|p(q_C)-p(q)|\le (1+C)^{-1}$,
\[
 \limsup_d\|\E W_C-p(Q_d)\|_{1,d}
 \le e_C:=|x_C-H_C(1)|+(1+C)^{-1}\longrightarrow0 .
\]
Restore all weights with $D^2=Q_dQ_{d,C}^{-1}$. The two multipliers in Lemma~\ref{lem:uncap} are positive contractions with disjoint supports in $q>C$ and $q<C^{-1}$. There $1-p(q)$ and $p(q)$, respectively, are at most $(1+C)^{-1}$. Their difference has norm at most one, so $\limsup_d\E b\le e_C+(1+C)^{-1}$. By \eqref{eq:spatial-tracenorm} and Jensen the limiting restoration cost is at most $2\sqrt{e_C+(1+C)^{-1}}$. Let $C\to\infty$.
At the endpoints use $\|U-V\|_{1,d}\le\min\{\tr_d(U+V),\tr_d(2I-U-V)\}$ for positive contractions. Monotonicity compares the traces at occupation $\epsilon$ when $x_d\to0$, and the deficits at $1-\epsilon$ when $x_d\to1$. If $\epsilon$ is unattainable, both traces are at most $\min(\rank Q_d,\rank K_d)/d\le\epsilon$. Now send $\epsilon\downarrow0$. This proves \eqref{eq:mean-matrix-track}
along arbitrary scales, including zero inputs.

For deterministic $0\preceq E_d\preceq I$, a row change costs at most
$2/d$ in normalized marked response by the graph estimate, since the matrix changes by rank at most two. Bounded differences give
summable tails about the means on a grid of at most $d+2$ scales of mesh
$1/d$ in $x_d(t)=d^{-1}\tr M_d(t)$, including endpoint limits.
Both marked traces increase, and each deterministic increment is
bounded by the full occupation increment. Interpolation gives uniform
almost-sure comparison, boundedness gives $L^1$, and positive and
negative parts handle symmetric contractions. The identity test and
\eqref{reb:two-scalar-accuracy} give \eqref{eq:track} and weak comparison
by Stieltjes uniqueness. Conversely either proposed comparison gives
the flat-test means used in Proposition~\ref{short:necessity}, by taking
traces or the identity test and using boundedness.
\end{proof}
\section{Singular values and recovery of the two inputs}
\subsection{Proportional singular values from the two input lists}\label{sec:rank-profile}
Order the spatial weights $q_i$ and time amplitudes $x_j=d\lambda_j$ decreasingly, padding either list by zeros. For $0<u<1$ put
\[
 E_d(u)^2=\frac2d\sum_{j>\lfloor du/2\rfloor}x_j^2,\qquad
 S_d(u)=q_{\lceil du\rceil} E_d(u),\qquad
 R_d(u)=s_{\lceil du\rceil}(A_d^X),
\]
where $s_i(A)$ denotes the $i$th singular value. The normalized temporal tail energy includes every remaining mode, even when the temporal rank is infinite.

\begin{theorem}[Proportional singular values]\label{thm:rank-profile}
Fix an input array satisfying \eqref{eq:pc}. Almost surely, for every $\delta>0$ and all sufficiently large $d$, whenever $0<a<b<1$ and $\Delta=b-a\ge\delta$,
\begin{equation}\label{eq:rank-profile}
 (\Delta^2/6)S_d(b)\le R_d(a),\qquad R_d(b)\le(6/\Delta^2)S_d(a).
\end{equation}
The ordered absolute roots satisfy the same bounds for all large $d$. For $U$ uniform on $(0,1)$, any joint subsequential limits of the empirical singular law and the law of $S_d(U)$ on $I_*$ have equal masses at zero and at infinity.
\end{theorem}
\begin{proof}
Write $T_d(u)=vy/u$ when $H_d(v)=G_d(y)=u$, and $T_d(u)=\infty$ when an input has normalized rank at most $u$. For $0<h<k<1$,
\begin{equation}\label{eq:profile-inverse-bounds}
 \frac{(k-h)^{3/2} S_d(k)}{\sqrt{hk}}\le T_d(h)^{-1},\qquad
 T_d(k)^{-1}\le\frac{k S_d(h)}{(k-h)^{3/2}}.
\end{equation}
At occupation $h$ the first $\lceil dk\rceil$ spatial terms give $1/v\ge (k-h)q_{\lceil dk\rceil}/h$. The first $j=1+\lfloor dk/2\rfloor$ temporal pairs give $1+y^2x_j^2\le k/(k-h)$; the tail starting at $j$ then gives $h\ge y^2E_d(k)^2(k-h)/k$. At occupation $k$, removing the terms before each rank-$h$ cutoff costs at most $h$ in each transform, so $vq_{\lceil dh\rceil}\ge k-h$ and $y^2 E_d(h)^2\ge k-h$. These bounds imply \eqref{eq:profile-inverse-bounds}, also for unavailable crossings by the rank convention.

Take $m=(a+b)/2$, $t=T_d(m)$. If $t=\infty$ then $R_d(b)=S_d(b)=0$. Otherwise the sample response $X(t)$ differs from $m$ by at most $\Delta/4$ eventually, uniformly for $\Delta\ge\delta$, by Theorem~\ref{thm:main}. Sorting bounds it between $bw_t(R_d(b))$ and $a+(tR_d(a))^2$. Thus $tR_d(a)\ge\sqrt\Delta/2$ and $tR_d(b)\le2/\sqrt\Delta$; \eqref{eq:profile-inverse-bounds} gives \eqref{eq:rank-profile} since $4\sqrt2<6$. For the roots use the error $4/d$. Along a weakly convergent subsequence, ranks strictly below the escaping mass diverge and higher ranks remain bounded; passing between profiles in \eqref{eq:rank-profile} gives equality of escaping masses. The threshold $1-\mu(\{0\})$ gives the corresponding zero-mass assertion.
\end{proof}
For any $L_d\to\infty$, $z_d\in\mathbb R$, the laws in $U$ of $(\log R_d(U)-z_d)/L_d$ and $(\log S_d(U)-z_d)/L_d$ approach each other weakly on $[-\infty,\infty]$ almost surely ($\log0=-\infty$). Indeed neighboring log quantiles on a fixed rank grid differ in the one-sided bounds by $O(1/L_d)$; compactify by $\arctan$ and refine the grid. This includes polynomial exponents $L_d=\log d$, $z_d=0$.

The gaps are necessary. At $u=1/2$, with Gaussian rows, $d=4n$, $Q=I$ and $d/4$ nonzero pairs give $R_d(u)>0=S_d(u)$ almost surely. Conversely $d=4n+2$, $d/2$ unit spatial weights (the rest zero) and $d$ unit pairs give $R_d(u)=0<S_d(u)$ by odd-dimensional skew symmetry.

\begin{corollary}[Trimming and the endpoints]\label{reb:two-boundaries}
Under \eqref{eq:pc} the root laws, the mean sample laws and almost surely the sample laws are tight on the finite half-line exactly when $S_d(u)$ stays bounded for each $0<u<1$. They converge to $\delta_0$ exactly when $S_d(u)\to0$ for every such $u$, and convergence in probability suffices for the samples. Both criteria are equivalent to the corresponding ones for $B_{d,m}=\sum_{i>m}q_i^2\sum_{j>m}\lambda_j^2$ at every $m=\lfloor du\rfloor$.
\end{corollary}
\begin{proof}
The rank bounds and Theorem~\ref{thm:main} give the endpoint tests for $S_d$. For $0<a<b/2$, summation up to or beyond the cutoffs gives $B_{d,\lfloor da\rfloor}\le S_d(a)^2/2$ and $S_d(b)^2\le 2 B_{d,\lfloor da\rfloor}/(b-a)$. Sufficiency even holds for general independent centered isotropic rows: deleting $m$ spatial directions and $m$ time pairs costs rank at most $4m$, leaving normalized mean squared Frobenius norm at most $2B_{d,m}$ by the degree-two coefficient. Apply Markov's inequality and rank-two bounded differences to bounded responses.
\end{proof}

\subsection{Recovery and spatial rank saturation}\label{sec:known-spatial}
Write $x_d=H_d(v_d)=G_d(y_d)$, $v_dy_d=tx_d$. In any reciprocal choice of units, if $x_d\to x>0$ and $v_d\to v>0$ on an interval, set $y=tx/v$. At one positive argument, $g_{y_d}\ge\min(1,y_d^2)$ bounds the source masses. Therefore every compactified input cluster $(\rho,\zeta)$ satisfies
\begin{equation}\label{eq:two-copy-transforms}
 H_\rho(v(t))=x(t)=G_\zeta(y(t)).
\end{equation}
Both inputs use the same compact transform $F_\nu(z)=\int(1+zu)^{-1}\,\nu(\dd u)$ on $[0,1]$. For $z>-1$,
\[
 F_{p_*\rho}(z)=1-\frac{zH_\rho(1+z)}{1+z},\qquad
 F_{(w_1)_*\zeta}(z)=\frac{G_\zeta(\sqrt{1+z})}{1+z}.
\]
Its Taylor coefficients at zero are the signed moments. Thus any interval of arguments determines each measure, endpoints included; varying $x$ and continuous $v$ provide such intervals for both.
Throughout the next theorem $\rho_d\Rightarrow\rho$ on the finite half-line and $r=1-\rho(\{0\})>0$.
\begin{theorem}[Recovery from one observed input]\label{thm:known-spatial}\label{thm:sat-bound}
Assume \eqref{eq:pc} and suppose $\eta_d\Rightarrow\eta+\alpha\delta_\infty$ on $I_*$
with $\eta$ supported on $[0,\infty)$.
\begin{enumerate}[label=\textup{(\roman*)}]
\item If $\alpha<r$, the output converges almost surely to a deterministic law $\mu$ with escaping mass $\alpha$. Its response $x(t)=\int w_t\,\dd\mu$ is zero for the zero source and otherwise satisfies
\[
 x(t)=H(v)=\alpha+G_\eta(tL(v)).
\]
The pair $(\rho,\mu)$ determines $(\eta,\alpha)$; if $\eta\ne0$, the pair $((\eta,\alpha),\mu)$ also determines $\rho$.
Without assuming source convergence, output convergence in probability to a deterministic law with mass at infinity below $r$ already implies convergence of the temporal sources.
\item If $\alpha=r$, then $\mu_d^{Q,X}\Rightarrow(1-r)\delta_0+r\delta_\infty$ almost surely, whatever the finite part $\eta$. If $\alpha>r$, every subsequential output limit almost surely has escaping mass $\beta\in[r,\min(1,\alpha)]$. For every prescribed limiting input pair $(\rho,\eta+\alpha\delta_\infty)$, every such $\beta$ is attained by some arrays with these limits.
\end{enumerate}
\end{theorem}
\begin{proof}
By Theorem~\ref{thm:main} it suffices to follow $x_d=H_d(v_d)=G_d(tL_d(v_d))$. At each $v>0$ the two sides tend to $H(v)$ and $\alpha+G_\eta(tL(v))$; $g_y$ is jointly continuous on compact positive $y$-intervals times $I_*$. For a nonzero source and $\alpha<r$ these monotone curves cross uniquely, so squeezing gives the stated response. For the zero source test arbitrarily small $v$. Stieltjes uniqueness gives weak convergence. Testing fixed $v$ with $H(v)>\alpha$ as $t\downarrow0$ gives $x(0+)=\alpha$.

Given $\rho$ and nonzero source, $v=H^{-1}(x)$ gives $y=tx/v$. Constant $y$ would force constant $x,v$ and hence contradict $y=tL(v)$. Thus \eqref{eq:two-copy-transforms} identifies the time source. Conversely, for known $\eta\ne0$, $y=G_\eta^{-1}(x-\alpha)$ gives $v=tx/y$; $v$ varies, recovering $\rho$. Without source convergence choose $x(t)<H(v)<r$ using $x(0+)<r$. Eventually $v_d<v$, so $G_d(tL_d(v))\le x_d$ at arguments tending to $tL(v)>0$. The source masses are bounded and each cluster has escaping mass at most $x(t)<r$; the preceding inverse identifies it.
For $\alpha\ge r$ the trial crossing forces $v_d\to\infty$, hence $\liminf x_d\ge r$ and $y_d=tx_d/v_d\to0$. Now $g_{y_d}\le1$ eventually and tends uniformly to zero on bounded sets, so $\limsup x_d\le\alpha$. At equality this identifies the two-atom limit. Otherwise send $t\downarrow0$ in each limiting response to bound its escaping mass.

To realize $r\le\beta\le\min(1,\alpha)$, approximate the positive part of $\rho$ on $\lfloor rd\rfloor$ directions, flooring at $d^{-1/2}$. Use the same floor up to rank $\lfloor\beta d\rfloor$, then zero. Take $\lfloor\alpha d/2\rfloor$ unit time pairs. Add pairs realizing $\eta$ by atomic approximation: each $m\delta_s$, $s>0$, takes $\lfloor dm/(2w_1(s))\rfloor$ pairs at $d\lambda=s$, and diagonalization includes $s=0$. Thus $S_d(u)\ge\sqrt{d(\alpha-u-o(1))}$ for $u<\beta$, whereas $S_d(u)=0$ for $u>\beta$. Theorem~\ref{thm:rank-profile} gives the claimed escape.
\end{proof}
Under Theorem~\ref{thm:known-spatial}, on any output-convergent subsequence the escaping mass equals $\sup\{u:q_{\lceil du\rceil}x_{\lceil du/2\rceil}\to\infty\}$ (empty supremum zero). Indeed $E_d(u)^2\le G_d(1)(1+x_{\lceil du/2\rceil}^2)$, the spatial quantiles stay bounded at fixed positive ranks, and $E_d(a)\ge\sqrt{b-a-2/d}\,x_{\lceil db/2\rceil}$ for $0<a<b<1$ and large $d$. Close the rank gaps in \eqref{eq:rank-profile}. For spatial thirds weighted $1,d^{-\kappa},0$ ($\kappa>0$) and $\lfloor d/2\rfloor$ unit time pairs the same bounds give $R_d(u)\asymp_u d^{1-\kappa}$ for $1/3<u<2/3$, with identical limiting inputs for all $\kappa$. At $Q=I$, $c=1-\alpha>0$, the normalized finite part has complementary response $q$ solving $q+G_{\eta/c}(ctq)=1$, since $L=1-H$.
\subsection{The spatial law from two independent samples}
Multiplying $Q_d$ by $a>0$ and dividing $K_d$ by $a$ leaves every sample unchanged. Take two independent samples with common inputs in the same spatial coordinates, and put
\[
 \widehat x_d(t)=\tr_d W_d(t),\qquad \widehat c_d(t,s)=\tr_d[W_d(t)W'_d(s)] .
\]
Only the second statistic uses the relative orientation of the two samples.

For $H=H_\rho$ put $C_\rho(v,w)=\int p(vq)p(wq)\,\dd\rho$. The elementary divided difference of $p$ gives
\[
 H(v)-C_\rho(v,w)=v\frac{H(v)-H(w)}{v-w},
\]
with diagonal value $vH'(v)$. When $\rho((0,\infty))>0$, use $u=H(v)\in(\rho(\{\infty\}),1-\rho(\{0\}))$ as coordinate and put
\begin{equation}\label{eq:spatial-derivative}
 A_\rho(u)=H(v)-C_\rho(v,v)=vH'(v),\qquad
 (\log H^{-1})'(u)=1/A_\rho(u).
\end{equation}
Equality of $A_\rho$ on any common open interval makes the inverse transforms proportional there. Analyticity and Stieltjes uniqueness identify $\rho$ up to dilation, which conversely preserves $A_\rho$.

\begin{theorem}[Spatial information in independent-sample overlaps]\label{thm:replica}
Suppose the rows of each sample satisfy the hypotheses of Theorem~\ref{thm:main}, including \eqref{eq:pc}. The two row laws may differ. Almost surely and in $L^1$,
\begin{equation}\label{eq:two-copy-forward}
 \sup_{t,s>0}\big|\widehat c_d(t,s)-\tr_d[N_d(t)N_d(s)]\big|\longrightarrow0 .
\end{equation}
Fix $s>0$. Suppose $\widehat x_d(t)$ and either the diagonal $\widehat c_d(t,t)$ or the section $\widehat c_d(t,s)$ converge in probability to deterministic limits on an open interval containing $s$, with $x=\lim\widehat x_d$ nonconstant. Then the spatial law $\rho$ and temporal source $\zeta$ of the pair normalized as $(a_dQ_d,K_d/a_d)$, $a_d=v_d(s)$, exist as weak limits on $I_*$ and are uniquely determined. The diagonal data $b(t)=c(t,t)$ give the first line below. Section data $c(t,s)$ give the quotient in the second. Both recover $r(t)=\lim v_d(t)/a_d$ and hence $y(t)=tx(t)/r(t)$.
\begin{align}
 x(t)-b(t)&=A_\rho(x(t)),&
 \log r(t)&=\int_s^t\frac{\dd x(z)}{x(z)-b(z)},\label{eq:diagonal-clock}\\
 r(t)&=\frac{x(t)-c(t,s)}{x(s)-c(t,s)},&
 y(t)&=\frac{t x(t)}{r(t)} .\label{eq:relative-parameter}
\end{align}
Observation on any set with an accumulation point in $(0,\infty)$ suffices if the limiting ordinary readings are nonconstant.
\end{theorem}
\begin{proof}
Independence factors the mean overlap, so \eqref{eq:mean-matrix-track} controls its bias. On the product grid of mesh $1/d$ in $\Theta_d$, each row edit costs $2/d$ and bounded differences give summable tails. Monotonicity in each argument and the marginal increment bounds give \eqref{eq:two-copy-forward} by interpolation. Replacing both $N_d$ by $M_d$ costs at most $8/d$ by \eqref{eq:finite-spatial-gap}.
Put $r_d(t)=v_d(t)/v_d(s)$, discarding finitely many zero inputs. For $t\ge s$, the parameters increase while $L_d$ and $G_d(y)/y^2$ decrease, giving $y_d(t)/y_d(s)\le t/s$ and $1\le r_d(t)\le(t/s)y_d(t)/y_d(s)\le(t/s)^2$.
Hence each subsequence admits a spatial cluster $\rho$ and locally uniform limit $r>0$, with $r(s)=1$ and $x(t)=H_\rho(r(t))$; the overlap in this cluster is $C_\rho(r(t),r(z))$. Nonconstant $x$ forces interior spatial mass. The divided difference above gives \eqref{eq:relative-parameter} by exchanging its arguments, and its diagonal gives \eqref{eq:diagonal-clock} by integration. Either reading fixes $r$ and \eqref{eq:two-copy-transforms} identifies both measures.

The probability measure on root pairs with rectangle masses $\tr_d[\Xi_d(J)\Xi_d(L)]$ has both marginals $\nu_d^Q$ and product response $\tr_d[N_d(t)N_d(z)]$. On $I_*^2$ its subsequential limits have analytic diagonal and fixed-section responses, since the kernels are locally uniformly bounded in a complex neighborhood of positive scales. The identity theorem extends these and the ordinary response from any accumulating set.
\end{proof}

For $r(t)\ne r(z)$ the same identity gives $c(t,z)=[r(t)x(z)-r(z)x(t)]/[r(t)-r(z)]$. At $u=x(t)$ the excess $b(t)-u^2=\Var_\rho[p(r(t)q)]$ vanishes at an interior occupation exactly for a point mass. On $[0,\infty)$ the identity $S_\rho(-u)=(1-u)H_\rho^{-1}(u)/u$ gives $(\log S_\rho(-u))'=(b(t)-u^2)/[u(1-u)(u-b(t))]$.

For $x\equiv a\in(0,1)$ the pair in anchor units is identifiable exactly when $b=c(s,s)=a^2$. Indeed $U=p(q)$ under any normalized spatial cluster has $\E U=a$, $\E U^2=b$. If $b<a$, strict increase of $H_\rho$ gives $r=1$, $G_\zeta(ta)=a$ and $\zeta=a\delta_\infty$. Any such law of $U$ occurs by empirical approximation of $p^{-1}(U)$ with $\lfloor ad/2\rfloor$ unit time pairs. At $b=a^2$ this fixes $q=a/(1-a)$; for $a^2<b<a$ the two-point laws on $\ell<a<u=a+(b-a^2)/(a-\ell)<1$ with mean $a$ give infinitely many solutions. At $b=a$ the spatial law is $(1-a)\delta_0+a\delta_\infty$. Both temporal sources $a\delta_j/g_{sa}(j)$, $j=0,1$, give $r=t/s$ and $c(t,s)=a$. Approximate the spatial atoms by $d^{-1},d$ and realize the sources as in Theorem~\ref{thm:sat-bound}.
\paragraph{Accuracy of recovery.}
Let $[s,T]$ be an observed interval, $s>0$. For two input pairs in anchor units, use their crossing readings $x_i(t)=H_{\rho_i}(r_i(t))$ and $c_i(t,z)=C_{\rho_i}(r_i(t),r_i(z))$, with $x_i(T)-x_i(s)\ge\delta>0$. Write $z_i(t)=c_i(t,t)$ for both models or $z_i(t)=c_i(t,s)$ for both. Define $e=\|x_1-x_2\|_\infty+\|z_1-z_2\|_\infty$ on $[s,T]$. Use bounded--Lipschitz distance on $[0,1]$ and put
\[
 D=d_{\rm BL}(p_*\rho_1,p_*\rho_2)+d_{\rm BL}((w_1)_*\zeta_1,(w_1)_*\zeta_2).
\]
These coordinates include the endpoint masses. With $C$ depending only on $s,T,\delta$,
\begin{equation}\label{eq:inverse-modulus}
 \|\log r_1-\log r_2\|_\infty\le Ce,\qquad
 D\le C/\log(2+1/e).
\end{equation}
At $e=0$ this means $D=0$. For exact finite root readings in dimensions $d_1,d_2$, use their discrepancy plus $12/d_1+12/d_2$, by \eqref{eq:finite-spatial-gap}. Proposition~\ref{prop:finite-overlaps} shows the logarithmic order is optimal even for errors measured on the entire overlap surface. The upper estimates are proved in Appendix~\ref{r21:sec-logarithms}.

\begin{corollary}[Absolute recovery from one calibrated probe]\label{cor:calibrated-inverse}
Assume \eqref{eq:pc}. For a deterministic real positive contraction $B_d$ suppose
$\sigma_d\Rightarrow\sigma$ on $I_*$, where $\sigma_d(J)=d^{-1}\tr B_d\mathbf1_J(Q_d)$, with known $\sigma$ charging $(0,\infty)$. Suppose the ordinary response and $d^{-1}\tr B_dW_d(t)$ converge in probability to deterministic $x(t),m(t)$ on an open interval. If $m$ takes two distinct values strictly between $\sigma(\{\infty\})$ and $\sigma((0,\infty])$, both input measures converge on $I_*$ and their limits are determined in the original units. On the interior range,
\[
 v(t)=H_\sigma^{-1}(m(t)),\qquad y(t)=tx(t)/v(t).
\]
For $\sigma=\gamma\delta_p$, $p>0$, this gives $v=m/[p(\gamma-m)]$.
The conclusions apply whenever $\rho_d\Rightarrow\rho$ on $[0,\infty)$ and $\eta_d\Rightarrow\eta+\alpha\delta_\infty$ on $I_*$ with $\eta\ne0$ and $\alpha<1-\rho(\{0\})$; then $\alpha=x(0+)$.
\end{corollary}
\begin{proof}
On an interior subinterval where $m$ varies, Theorem~\ref{thm:main} gives $H_{\sigma_d}(v_d)\to m$. Squeezing against the strictly increasing $H_\sigma$ gives the displayed continuous $v$. Every spatial cluster dominates $\sigma$, so $x=H_\rho(v)$ varies and $x\ge m>0$. Apply \eqref{eq:two-copy-transforms}. In the stated subcritical case Theorem~\ref{thm:known-spatial} gives strictly increasing $x,v$, hence $m=H_\sigma(v)$ by Theorem~\ref{thm:main}.
\end{proof}
\begin{corollary}[One observation for a fixed kernel]\label{thm:fixed}
Let $Q_d=I_d$, $K_d=K$, and assume \eqref{eq:pc}. At any prescribed $t_0>0$, convergence in probability of the sample response to a deterministic $x_0<1$ is equivalent to $j\lambda_j(K)\to c$ for some finite $c\ge0$. In that case the full empirical law converges almost surely to the law of $\pi c|C|$, $C$ standard Cauchy, and $x_0/(1-x_0)=\pi c t_0$. In particular these conditions are equivalent to convergence in probability to a deterministic probability law on the finite half-line.
\end{corollary}
The proof is combined with the endpoint assertions in Appendix~\ref{r21:sec-dilation}.
\section{Brownian area and its time mesh}\label{r21:applications}\label{reb:area}
Put $b=T/2>0$. Brownian area $\tfrac12\int_0^T(W^i\,\dd W^j-W^j\,\dd W^i)$ with covariance $Q_d$ uses the Gaussian rows and kernel $S(s,u)=\tfrac12\operatorname{sgn}(u-s)$. Compressing to normalized indicators of cells $h_j=Tp_j$ gives polygonal area, with row coordinates $\Delta_j B^i/\sqrt{h_j}$ and $K_{jk}=b\sqrt{p_jp_k}\operatorname{sgn}(k-j)$. Both sample laws below allow independent rows satisfying \eqref{eq:pc}.

\begin{corollary}[Spatial observations of Brownian area]\label{reb:brownian}
For standard Cauchy $C$, set $\HC_a=\mathcal L(a|C|)$, including point masses at $a=0,\infty$, and $\mathcal B_b(\rho)=\int\HC_{bq}\,\rho(\dd q)$ for finite measures. For $S$ and arbitrary moving $Q_d$,
\begin{equation}\label{reb:brownian-quantile}
 d_{\rm K}(\nu_d^Q,\mathcal B_b(\rho_d))\le 1/d,\qquad
 d_{\rm K}(\mu_d^{Q,X},\mathcal B_b(\rho_d))\to0\quad\text{almost surely}.
\end{equation}
The first bound is attained for every $Q_d\ne0$. For any deterministic real $0\preceq B_d\preceq I$ the second comparison also holds for the weighted laws $\mu_{d,B}(J)=\tr_d[B_d\mathbf1_J(|A_d^X|)]$, $\rho_{d,B}(J)=\tr_d[B_d\mathbf1_J(Q_d)]$. The mixture map is a homeomorphism onto a closed subset of $\mathcal P(I_*)$. Thus the spatial laws converge on $I_*$ exactly when the samples converge almost surely to a deterministic law, or when their means converge.
\end{corollary}
\begin{proof}
The equations $\zeta f'=-f$, $f(T)=-f(0)$ give $\lambda_j=T/((2j-1)\pi)$. Hence $\prod_j(1+s^2\lambda_j^2)=\cosh(bs)$, $(2k)!e_k(\lambda^2)=b^{2k}$ and $p_d(x)=\operatorname{Re}\prod_i(x+ibq_i)$ \cite{Zhao}. Its positive roots mark the half-integer crossings of the decreasing phase $\phi/\pi=\pi^{-1}\sum_i\arctan(bq_i/x)$. Counting each twice rounds the mixture tail $2\phi/(\pi d)$ within $1/d$, attained at a positive root or at zero for rank one.

Here $P_d(t)$ is the even part of $\prod_i(1+btq_i)$. Thus the eigenvalues of $N_d(t)$ are the inclusion probabilities of independent Bernoulli variables of means $p(btq_i)$ conditional on even total. The even and odd laws couple at Hamming distance one: couple the last bits, using identical prefixes when they differ and an inductive opposite-parity coupling when they agree. Deterministic bits may be omitted. Mixing the two parities gives $\|N_d(t)-p(btQ_d)\|_{\mathcal S_1}\le\Pr(\mathrm{odd})\le1$. Since $\int w_t\,\dd\HC_{bq}=p(btq)$, Theorem~\ref{thm:main} transfers all marked responses. A Kolmogorov discrepancy would persist at one after dilation by discrepant positive thresholds (right continuity covers zero). Paired compactified limits agree by their responses and masses, and have no atom at one by the uniform log-density bound $1/\pi$ for positive Cauchy components. This proves the comparison, also for means by bounded convergence. The mixture map is weakly continuous since $0<|C|<\infty$ almost surely; Stieltjes uniqueness for its complementary response $\int(1+btq)^{-1}\rho(\dd q)$ and compactness give the homeomorphism.
\end{proof}

\begin{corollary}[Time-weighted mesh recovery]\label{prop:joint-mesh}
Let $\nu_d=\sum_jp_j\delta_{dp_j}$ be the law of the rescaled cell size seen at a uniform time. Suppose $\rho_d\Rightarrow\rho$ on $[0,\infty)$ with $\rho((0,\infty))>0$. If $\nu_d\Rightarrow\nu$ on $I_*$, the polygonal sample laws converge almost surely to a deterministic law on $[0,\infty)$ with response
\begin{equation}\label{eq:joint-mesh-scales}
 t=\frac{H_\rho^{-1}(y)B_\nu^{-1}(y)}{y},\qquad
 H_\rho(v)=\int\frac{vq}{1+vq}\,\rho(\dd q),\qquad
 B_\nu(w)=bw\int\frac{\nu(\dd u)}{1+bwu}.
\end{equation}
Here $0<y<\min(H_\rho(\infty),B_\nu(\infty))$ for $\nu\ne\delta_\infty$, and $y=0$ otherwise. Conversely convergence in probability of the output to a deterministic law implies profile convergence. Given $\rho$ or $\nu\ne\delta_\infty$, the output identifies the other input. At any fixed $t>0$,
\[
 \int w_t\,\dd\mu_d^{Q,X}\xrightarrow{\PP}H_\rho(bt)
 \quad\Longleftrightarrow\quad
 \sum_{j:\,d h_j>\varepsilon}h_j\to0\quad(\varepsilon>0).
\]
\end{corollary}
\begin{proof}
First-row expansion gives Pfaffian one for each even principal minor of the sign matrix, hence mean coefficient $(2k)!b^{2k}e_{2k}(q)e_{2k}(p)$. At the crossing for $\rho_d,\nu_d$, Bernoulli counts with odds $vq_i,dwbp_j$ have common mean $x=dy$. Their masses $s_n,h_n$ give degree weights proportional to $(2k)!s_{2k}h_{2k}/x^{2k}$; Theorem~\ref{reb:mode}'s stride-two comparison yields $|d\Theta_d-x|<4$.

For $\nu\ne\delta_\infty$, $H_\rho^{-1}(y)/y$ is nondecreasing, both inverses are $O(y)$ at zero, and one diverges at the upper endpoint. Thus the inverse scale runs from zero to infinity and convergence of crossings follows by squeezing. For $\nu=\delta_\infty$, $B_{\nu_d}(w)\to0$ gives zero limiting occupation. Theorem~\ref{thm:main} transfers the limit; $y(0+)=0$ excludes escape. For a nonzero response either input fixes the other transform on an interval, and $B_\nu(1/(bv))=\int(v+u)^{-1}\,\dd\nu$ and $1-H_\rho(1/v)=v\int(v+q)^{-1}\,\dd\rho$. Stieltjes uniqueness and compactness give the converse. Only $\nu=\delta_0$ attains $B_\nu(w)=bw$ at any $w>0$, equivalently $y(t)=H_\rho(bt)$ at one scale. This is the stated cell-size criterion.
\end{proof}

For $s=\sum_iq_i>0$, exchanging $(q,p)$ with $(sp,q/s)$ preserves the mean polynomial in padded dimension $\max(d,N)$ for $N$ cells. Take $N=d\in4\mathbb N$, $Q=I$, and $dp_j=1/2,1,3/2$ in proportions $1/4,1/2,1/4$. The calibrated middle block ($\tfrac12\delta_1$) gives $x(t)/2$ before exchange and $p(v)/2>x(t)/2$ afterwards by strict concavity. Their same-scale overlaps also differ by the spatial variance criterion. Appendix~\ref{app:mesh} identifies $A_\rho$ as the leading mesh-loss coefficient.
\begin{proposition}[Resolution of the two inputs]\label{prop:finite-overlaps}
Given readings at any $k$ positive scales, there are two different normalized input pairs with the same ordinary curve and identical limits for all overlaps between these scales. In mean-one units the spatial weights and rescaled cell sizes may each take just $2k+3$ possible limiting values in $[1/2,3/2]$. Even with all scales observed, the order $1/\log(1/e)$ in \eqref{eq:inverse-modulus} cannot be improved uniformly.
\end{proposition}
\begin{proof}
Take $b=1$, $N=d$ and let $\rho,\xi$ be the limiting number-weighted laws of $q_i,dp_j$, both of mean one. Then $\nu(\dd u)=u\,\xi(\dd u)$ and \eqref{eq:joint-mesh-scales} gives $t=H_\rho^{-1}(x)H_\xi^{-1}(x)/x$, invariant under exchange. For $n\ge2k+2$, use $q_j=1/2+j/n$ with baseline masses $\omega_j=2^{-n}\binom nj$, $0\le j\le n$. Let $v_i$ be the spatial parameters at the prescribed scales when both inputs have this baseline law. Define
\[
 P(q)=\prod_{i=1}^k\left(\frac{1+v_iq}{1+3v_i/2}\right)^2,\qquad
 \rho_\pm(\{q_j\})=\omega_j[1\pm\tfrac12(-1)^jP(q_j)].
\]
Alternating binomial weights annihilate polynomials of degree $<n$. Multiplication by $P$ makes $1,q,p(v_iq)$ and the tested products such polynomials. Thus $\rho_\pm$ are mean-one probabilities; exchanging them preserves the entire ordinary curve and the prescribed overlaps. Equal means exclude a nonunit dilation. Choose empirical frequencies and rescale the weights to mean one.

Uniformly in $v,w>0$, $P(q)p(vq)p(wq)$ is bounded by $C_k$ on $|q-1|\le2/3$. Its degree-$(n-1)$ Taylor error on the support is $O_k((3/4)^n)$, bounding the overlap discrepancy at common parameters. For $H_+-H_-$ this gains $\min(v,1/v)$ by using $Pp(vq)$ for $v\le1$ and $P(1-p(vq))$ otherwise. At a common occupation the inverse parameters differ by at most a factor three by the support bounds. Since $\partial_{\log v}H_\pm(v)\asymp\min(v,1/v)$, their log difference is $O_k((3/4)^n)$. Lipschitz continuity of $p(vq)$ in $\log v$ gives the same bound on the entire overlap surface.

Interpolate $(-1)^j/(2n)$ at $q_j$ to a $1$-Lipschitz test. Its integrals differ by at least $3^{-2k}/(2n)$. Both laws tend to $\delta_1$ by binomial concentration, hence their anchor scales have a common positive finite limit and differ by $O_{k,s}((3/4)^n)$. Since $p$ is bi-Lipschitz on the resulting compact positive support, normalized spatial distance is at least $c_{k,s}/n$. The limiting curve $t=x/(1-x)^2$ gives positive oscillation on every fixed $[s,T]$. Take $n=2k+2$ for finite ambiguity and $n\to\infty$ for sharpness.
\end{proof}
\appendix
\section{The isotropic source and its inverse}\label{r21:canonical-section}
Take $Q_d=I$, $\mu_d=\mu_d^{I,X}$ and $\widehat q_d(z)=\int(1-w_z)\,\dd\mu_d$. For a finite positive measure $\eta$ on $[0,\infty)$ let $q_\eta(z)\in(0,1]$ solve $q+G_\eta(zq)=1$, and write $q_d=q_{\eta_d}$.

\begin{proposition}[The canonical law]\label{r21:law}
The map $\eta\mapsto\mu_\eta$ specified by
\begin{equation}\label{r21:canonical-law}
 \int(1-w_z)\,\dd\mu_\eta=q_\eta(z)\qquad(z>0)
\end{equation}
is a homeomorphism onto the weakly closed class of folded symmetric freely infinitely divisible laws.
\end{proposition}
\begin{proof}
Free L\'evy--Khintchine \cite[Theorem~3.2]{BP}, \cite[Theorem~1.1]{BG} parametrizes symmetric freely infinitely divisible laws $\nu_\eta$ by the symmetric extensions of $\eta$. The reciprocal Cauchy transform satisfies on $\mathbb C_+$
\begin{equation}\label{r22:free-canonical-transform}
 F+\varphi\circ F=\mathrm{id},\qquad \varphi(w)=\int\frac{w(1+s^2)}{w^2-s^2}\,\eta(\dd s).
\end{equation}
For the folded law, putting $F(i/z)=i/(zu)$ gives $u=\int(1-w_z)\,\dd\mu=q_\eta(z)$. Recovery follows as in Theorem~\ref{thm:known-spatial}, using
\begin{equation}\label{r21:stieltjes-source}
 G_\eta(v^{-1/2})=\eta([0,\infty))+(1-v)\int(v+s^2)^{-1}\,\eta(\dd s).
\end{equation}
Weak source limits transfer by monotone squeezing and $q_\eta(0+)=1$. Conversely a finite-line output limit has $q(0+)=1$, so $\limsup_d G_d(z/2)\le1-q(z)$ for small $z$. Source tightness follows from $g_z\ge\min(1,z^2)$ and
\begin{equation}\label{r21:source-tightness}
 \eta_d((M,\infty))\le2G_d(1/M),\qquad G_d(t)\le t^2(1+M^2)G_d(1)+\eta_d((M,\infty))
\end{equation}
for $M\ge1$, $t\le1$. Injectivity identifies every cluster point.
\end{proof}

\begin{corollary}[Moving spectra and escape]\label{r21:moving}\label{r25:escape}\label{r21:compact}
Under \eqref{eq:pc},
\begin{equation}\label{r21:moving-equivalent}
 \sup_{z>0}|\widehat q_d(z)-q_d(z)|\longrightarrow0\qquad\text{almost surely and in }L^1.
\end{equation}
Consequently $\mu_d$ follows $\mu_{\eta_d}$ almost surely in any weak metric on $I_*$. If $\eta_d\Rightarrow\eta+\alpha\delta_\infty$ there and $c=1-\alpha>0$, then
\[
 \mu_d\Rightarrow\alpha\delta_\infty+c(D_c)_*\mu_{\eta/c}\quad\text{almost surely},\qquad D_c(s)=cs.
\]
A deterministic output limit in probability with positive finite-line mass determines a unique compactified source limit. Complete escape in probability is equivalent to $\liminf_dG_d(t)\ge1$ for every $t>0$ and is then almost sure. This includes $\alpha\ge1$ or $G_d(1)\to\infty$.

On the finite line, weak convergence of $\eta_d$, pointwise convergence $G_d\to G$ with $G$ finite and $G(0+)=0$, weak convergence of $\E\mu_d$ to a probability, and weak convergence in probability of $\mu_d$ to a deterministic probability are equivalent. The limits are $\eta,G_\eta,\mu_\eta,\mu_\eta$. Transform convergence is locally uniform on $[0,\infty)$. Tight mean laws suffice for $\lim_{t\downarrow0}\limsup_d G_d(t)=0$ and $d_{\mathrm L}(\mu_d,\mu_{\eta_d})\to0$ almost surely.
\end{corollary}
\begin{proof}
Apply Theorems~\ref{thm:main} and~\ref{thm:known-spatial} at $Q=I$. For $c>0$ the limit solves $q+G_\eta(zq)=c$, giving $q=cq_{\eta/c}(cz)$ and the folded free compression $D_c(\nu_\eta^{\boxplus1/c})$ by \eqref{r22:free-canonical-transform} \cite[Sections~1.11--1.14]{NS}. Testing positive trial $q$ gives complete escape under the stated $\liminf$; conversely $G_d(t)\le1-\epsilon$ implies $q_d(t)\ge\epsilon$. The tail bound $\mu([1/z,\infty))\le2\int w_z\,\dd\mu$ transfers tightness, after which Proposition~\ref{r21:law} identifies paired limits. For transform convergence, \eqref{r21:source-tightness} and $G(0+)=0$ give tightness, and monotonicity gives local uniformity.
\end{proof}
With $d^2$ pairs of size $\sigma/(\sqrt2\,d^{3/2})$, $\eta=\sigma^2\delta_0$ gives a folded semicircle of variance $\sigma^2$. Adding $\lfloor\alpha d/2\rfloor$ unit pairs, $0<\alpha<1$, changes its normalized finite-part variance parameter to $(1-\alpha)\sigma^2$.
\begin{corollary}[Two-array recovery]\label{r23:two-array-recovery}
For two arrays under \eqref{eq:pc} with the first mean laws tight, under any joint coupling,
\begin{equation}\label{r21:two-arrays}
 d_{\mathrm{BL}}(\eta_d^K,\eta_d^L)\to0
 \iff \sup_{0\le t\le T}|G_d^K(t)-G_d^L(t)|\to0\ (T>0)
 \iff d_{\mathrm L}(\mu_d^K,\mu_d^L)\xrightarrow{\PP}0 .
\end{equation}
Here $d_{\mathrm{BL}},d_{\mathrm L}$ are the bounded-Lipschitz and L\'evy metrics. Each condition gives tightness of the second array and almost-sure output comparison.
\end{corollary}
\begin{proof}
Tightness transfers by Lipschitz cutoffs (including the constant one), \eqref{r21:source-tightness}, or the L\'evy inequalities. On tight pairs Proposition~\ref{r21:law} and Corollary~\ref{r21:moving} identify the cluster points.
\end{proof}
For an isometry $U:\mathbb R^m\to\mathbb R^d$ with $m\asymp d$ and $Q=UU^*$, $e_j(Q)=\binom mj$ factors $x^{d-m}$ out of the mean polynomial. Under \eqref{eq:pc}, removing these zero modes shows that $U^*X_dK_dX_d^{\mathsf T}U$ tracks the $m$-dimensional root law: the normalized response error is multiplied by $d/m$.
\section{Logarithmic observables and finite-dimensional bounds}\label{r21:sec-logarithms}
The logarithmic estimates use $Q=I$ unless a general covariance is specified. Keep the notation of Appendix~\ref{r21:canonical-section}. Write
$F_d(z)=d^{-1}\log\det(I+zA_d^X)$ and $x_k=d\lambda_{k,d}$. Set
\begin{align}
 V_\eta(z)&=\min_{q>0}\{q-\log q-1+\Phi_\eta(zq)\},\qquad
 q_\eta+G_\eta(zq_\eta)=1,\label{r21:root}\\
 \Phi_\eta(t)&=\int\frac{1+x^2}{2x^2}\log(1+t^2x^2)\,\eta(\dd x),
 \qquad \Phi_d(t)=d^{-1}\sum_k\log(1+t^2x_k^2).\label{r21:V-intro}
\end{align}
The integrand equals $t^2/2$ at zero; write $V_d=V_{\eta_d}$ and $\Phi_d=\Phi_{\eta_d}$.
If $\Phi_\eta(1)<\infty$, then $t\Phi'_\eta(t)=G_\eta(t)$, so $q_\eta$ is the unique minimizer. Differentiating at $t=zq_\eta(z)$, $B=tG'_\eta(t)$, $u=\log z$, and testing the minimum at one gives
\begin{equation}\label{r21:response-calculus}
 \partial_u V_\eta=1-q_\eta,\qquad
 -\partial_u q_\eta=q_\eta B/(q_\eta+B)\le1,\qquad
 e^{-1-\Phi_\eta(z)}\le q_\eta(z)\le1.
\end{equation}
\begin{proposition}[The scalar potential]\label{r21:saddle}
Whenever $\Phi_d(z)\le B$,
\begin{equation}\label{r21:saddle-est}
 \left|d^{-1}\log P_d(z)-V_d(z)\right|\le C_B/d .
\end{equation}
\end{proposition}
\begin{proof}
Put $w=dq_d(z)$ and $D(y)=\prod_k(1+y^2\lambda_k^2)$. Let $N\sim\mathrm{Pois}(w)$ and independently $S=2\sum_k B_k$ for independent Bernoulli variables of odds $(zw\lambda_k)^2$. Coefficient extraction in $e^yD(zy)$ gives
$P_d(z)=d!w^{-d}e^wD(zw)\PP(N+S=d)$.
At the crossing $\E(N+S)=d$, $\Var S\le2d$ and $e^{-1-B}d\le w\le d$. With probability at least $7/8$, $|d-S-w|\le4\sqrt d$; throughout this window the Poisson mass is comparable to $d^{-1/2}$ for large $d$, also its order at the mode. For bounded $d$ use $N=d,S=0$ with $\PP(S=0)\ge e^{-dB}$. Stirling's formula gives \eqref{r21:saddle-est}.
\end{proof}

By \eqref{r21:source-tightness}, boundedness and tightness of the sources are equivalent to the first of the following conditions:
\begin{align}
 \lim_{t\downarrow0}\limsup_d G_d(t)&=0,\tag{RC}\label{r21:rc}\\
 \lim_{M\to\infty}\limsup_d d^{-1}\sum_{x_k>M}\log(1+x_k^2)&=0.\tag{LI}\label{r21:li}
\end{align}

\begin{theorem}[Logarithmic extension]\label{r21:logarithms}
Under \eqref{eq:pc}, \eqref{r21:rc} and \eqref{r21:li}, $F_d-V_d\to0$ almost surely and in $L^1$, locally uniformly in $\log z$. The sources converge weakly to a finite measure if and only if $\Phi_d(t)$ converges to a finite number at each $t>0$. Along a subsequence with $\eta_d\Rightarrow\eta$, $\Phi_d\to\Phi_\eta$ locally uniformly on $[0,\infty)$ and $V_d(z),F_d(z)$ converge to
\begin{equation}\label{r21:log-law}
 V_\eta(z)=\tfrac12\int\log(1+z^2s^2)\,\mu_\eta(\dd s),
\end{equation}
where $\mu_\eta$ is the isotropic output with response $1-q_\eta$.
\end{theorem}
\begin{proof}
Since $(d)_{2j}\le d^{2j}$, $\E e^{dF_d(t)}=P_d(t)\le e^{d\Phi_d(t)}$. Thus $\E F_d(t)\le\Phi_d(t)$ and $\limsup_d(F_d(t)-\Phi_d(t))\le0$ almost surely by exponential Markov bounds. Also $0\le V_d(t)\le\Phi_d(t)$. Splitting the source integral at $M$ gives $\lim_{t\downarrow0}\limsup_d\Phi_d(t)=0$ by \eqref{r21:rc}--\eqref{r21:li}. On $[t,z]$ the integrated derivative error is bounded by $|\log(z/t)|\sup_s|\widehat q_d(s)-q_d(s)|$. Apply \eqref{r21:moving-equivalent} and send $t\downarrow0$ through a countable sequence.
Under weak source convergence, \eqref{r21:li} permits uniform truncation in \eqref{r21:V-intro} for bounded $t$, giving convergence of the potentials and minima. Conversely \eqref{r21:rc} gives compactness, and each cluster potential determines $G=t\Phi'$ and hence its source by \eqref{r21:stieltjes-source}. Integrating \eqref{r21:canonical-law} from zero proves \eqref{r21:log-law}.
\end{proof}
A single Gaussian pair of size $e^{\beta d}$, $\beta>0$, has bulk $\delta_0$ but $F_d(z)\to2\beta$ almost surely: its determinant is $1+z^2e^{2\beta d}Y_d$, where orthogonal projection gives $Y_d$ the product of independent $\chi^2_d,\chi^2_{d-1}$ variables and hence $Y_d/d^2\to1$ by exponential tails.

For rates assume
\begin{equation}\tag{QC4}\label{eq:qc}
 \Var(X_i^*T X_i)\le\kappa\|T\|_{\HS}^2\qquad(T=T^*\text{ finite rank}).
\end{equation}
Then $\Delta_d(a)\le\sqrt{\kappa/(ad)}$. Examples include independent standardized coordinates of bounded fourth moment, their orthogonal images, and uniform rows on the radius-$\sqrt p$ sphere in $\mathbb R^p$ ($\kappa=2$). Rank-one testing gives the small-ball constant $1/[8(1+\kappa)]$ by Paley--Zygmund.

\begin{proposition}[Bias, variance and window error]\label{r21:schur}\label{r21:window}
For a real skew kernel $T$, write $S=T^*T$, $b=d(\tr S)^2+d^2\tr(S^2)$, and form $F_d,P_d,V_d$ using $T$. Under \eqref{eq:qc},
\begin{equation}\label{eq:direct-jensen}
 0\le d^{-1}\log P_d(z)-\E F_d(z)\le C_\kappa z^4 b,
 \qquad \Var F_d(z)\le C_\kappa z^4 b/d .
\end{equation}
In particular, when $d\|T\|_{\HS}^2\le E$ and $d\|T\|_{\op}\le M$,
\begin{equation}\label{r21:window-est}
 \E|F_d(z)-V_d(z)|\le C/d,\qquad
 \E|\widehat q_d(z)-q_d(z)|\le C/\sqrt d ,
\end{equation}
with $C$ depending on $E,M,\kappa,z$, locally uniformly in $\log z$.
\end{proposition}
\begin{proof}
For $A_J=\sum_{i\in J}(TX_i)(TX_i)^*$, $m=|J|$, expansion of the Gram entries gives
$\E\tr A_J^2\le m(\tr S)^2+m(\kappa+m-1)\tr S^2$.
The row update in Proposition~\ref{thm:q-transfer} gives $M_i\preceq z^2A_{\ne i}$. For a conditional mean the coefficient $C_j$ is its average weighted by the determinant with row $j$ deleted. Future rows may be integrated as independent signs in a temporal eigenbasis since the determinant is quadratic in each row. Each tilted future row has conditional density $(1+g^*Dg)/(1+\tr D)$, $D\succeq0$, by the same update. Its second moment is at most $3I$ by Cauchy--Schwarz and $\E(v^*g)^4\le3\|v\|^4$. Hence $C_j\preceq z^2(A_{<j}+3(d-j)TT^*)$. Thus both have expected squared Hilbert--Schmidt norm at most $C_\kappa z^4 b$. The conditional Jensen gap of $\log(1+X_j^*C_jX_j)$ is at most $\kappa\tr C_j^2/2$ by the second derivative bound, and the conditional variance for $M_i$ is at most $\kappa\tr M_i^2$ by Lipschitz continuity. Sum the gaps and tensorize the variances to get \eqref{eq:direct-jensen}; the $L^2$ truncation in Proposition~\ref{thm:q-transfer} covers infinite rank.
In the window $b=O(d^{-1})$ and $\Phi_d(z)\le z^2E/2$, giving the log bound by \eqref{r21:saddle-est}. Taking traces in \eqref{eq:finite-spatial-gap} bounds the response bias by $\sqrt{J_d(z)}+4/d$; bounded differences with row cost $2/d$ control its centered error.
\end{proof}

For $B_0=\sup_d G_d(1)<\infty$, deleting $x_k>M\ge1$ costs at most $(1+M^{-2})\eta_d((M,\infty))$ in each response, by rank and by the unit lower bound on the crossing slope. The retained energy is at most $(1+M^2)B_0$, giving
\begin{equation}\label{r21:cutoff-key}
 \E|\widehat q_d(z)-q_d(z)|\le2(1+M^{-2})\eta_d((M,\infty))+C_{B_0,M,\kappa,z}/\sqrt d .
\end{equation}
On a band $0<\varepsilon\le d\lambda\le M$, put $B=idT$, $C=Z^*Z/d$ for the restricted row matrix. Taking Schur complements in the block matrix with diagonal $(B^{-1},tI)$ and upper block $Z^*/\sqrt d$ gives $\#\{j:\lambda_j(iZTZ^*)>t\}=n_-(B^{-1}-C/t)-n_-(B^{-1})$ for $t>0$.

\begin{proposition}[Weak-trace rates]\label{r21:weak-rates}
If $A=\sup_{d,k}k\lambda_{k,d}<\infty$ and \eqref{eq:pc} holds, the conclusions $F_d-V_d\to0$ almost surely and in $L^1$ hold locally uniformly in $\log z$. The centered error satisfies
\begin{equation}\label{r21:weak-selfaverage}
\E|F_d(z)-\E F_d(z)|\le C_{A,z}\log(d+2)/\sqrt d.
\end{equation} Under \eqref{eq:qc},
\begin{equation}\label{r21:rate-log}
 \E|F_d(z)-V_d(z)|\le C_{A,\kappa,z}d^{-1/4}\log^{3/4}(d+2)
\end{equation}
with locally uniform constants.
\end{proposition}
\begin{proof}
The envelope $d\lambda_{\lceil dx\rceil,d}\le A/x$ gives $G_d(t),\Phi_d(t)\le\pi At$ and bounds the \eqref{r21:li} tail by $\int_0^{A/M}\log(1+A^2/x^2)\,\dd x$. Apply Theorem~\ref{r21:logarithms}; the variance bound in Proposition~\ref{short:necessity} gives \eqref{r21:weak-selfaverage} since $\|K\|_{\HS}^2\le2A^2\sum_{k\ge1}k^{-2}$.
Remove the first $m$ pairs, writing $A_d^X=U+B$ with $B$ the removed part. Factoring the determinant ratio and its reciprocal with the skew resolvents (norm at most one) gives
\[
 d|F_{U+B}-F_U|\le\sum_{j\le2m}\log(1+zs_j(B))\le\log\det(I+zB)+m\log2 .
\]
The last inequality uses $(1+x)^2\le2(1+x^2)$. Thus the mean change is at most $\Phi_d^{\rm head}(z)+(m/d)\log2$ and the change in $V_d$ at most $\Phi_d^{\rm head}(z)$ by testing the truncated minimizer. The envelope and $m!\ge(m/e)^m$ bound $\Phi_d^{\rm head}(z)$ by $O_{A,z}((m/d)\log(e+d/m))$. For the residual kernel, $\|T\|_{\HS}^2\le CA^2/m$, $\tr((T^*T)^2)\le CA^4/m^3$; \eqref{eq:direct-jensen} and \eqref{r21:saddle-est} give bias $O_{A,\kappa,z}(d/m^2+d^2/m^3+1/d)$. Take $m=\lceil d^{3/4}/\log^{1/4}(d+2)\rceil$.
\end{proof}

\subsection{Gaussian comparison and rate obstructions}
The square-root dependence for the determinant tilt in \eqref{eq:finite-spatial-gap} is sharp at fixed dimension. Take $d=3$, $H=2^{-1/2}\left(\begin{smallmatrix}1&1\\1&-1\end{smallmatrix}\right)$, $U=H\otimes H$ and $K=\left(\begin{smallmatrix}0&1\\-1&0\end{smallmatrix}\right)\otimes I_2$. Choose the rows independently as twice a uniformly random signed column from the respective bases $I_4,U,\operatorname{diag}(1,1,1,-1)U$. They are centered and isotropic. Between different bases all entries of $K$ have modulus $1/2$. At $Q=I$ all three $|A_{ij}|$ for $i<j$ equal two, so $D=1+12t^2$ is constant while $W_{12}=\pm w_t(\sqrt{12})/3$ with equal probabilities. Vary $Q=e^{hB}$ with $B=E_{12}+E_{21}$. Writing $J(h)$ for the gap and $V=\tr BW$, Taylor expansion of $\log\E D_h-\E\log D_h$ gives $J(h)=h^2\Var(V)/(2d)+O(h^3)$ and $J'(h)=h\Var(V)/d+O(h^2)$. Since $dJ'(h)=\tr B(N_d-\E W_d)$ at that input, the difference is at least a positive multiple of $\sqrt{J(h)}$.

\begin{proposition}[Gaussian rate for arbitrary inputs]\label{prop:rate}
For Gaussian rows and every $d\ge2$, with probability at least $1-2d^{-3}$,
\begin{equation}\label{eq:two-spectrum-tracking}
 \sup_{t>0}\left|\int w_t\,\dd\mu_d^Q-\Theta_d(t)\right|
 \le\min\{1,\sqrt{20\log(d)/d},35d^{-1}\log(d)(1+\log d)\}.
\end{equation}
\end{proposition}
\begin{proof}
The coefficients $B_k=[t^{2k}]\det(I+tA_d)\ge0$ have means $a_{k,d}$. For positive means the exterior small-ball bound \cite[equation~(14)]{Zhao} gives $\PP(B_k\le ua_{k,d})\le8k\sqrt{u/\pi}$ ($k\ge1$). Markov and a union bound put all positive ratios $B_k/a_{k,d}$ in $[d^{-10},d^{10}]$ with the stated probability. Hence $|h|\le\varepsilon=10\log(d)/d$ on $\mathbb R$, where $h(u)=d^{-1}\log[P_d(e^u)/\det(I+e^uA_d)]$ and $h'$ is the response difference. The two log curves are convex with second derivatives at most $1/2$, so secants of width $\sqrt{8\varepsilon}$ yield $|h'|\le\sqrt{2\varepsilon}$.
Paired roots give analyticity on $|\operatorname{Im}u|<\pi/2$. On $|\operatorname{Im}u|\le a=\pi/4$, $|w_{e^u}(s)|\le1$, hence $|h'|\le2$ and $|h|\le B=\varepsilon+2a$. Three-lines interpolation bounds $|h(x+iy)|$ by $\varepsilon(B/\varepsilon)^{|y|/a}$. Cauchy's estimate at radius $a/(1+\log(B/\varepsilon))$ gives $\|h'\|_\infty\le(e/a)\varepsilon[1+\log(B/\varepsilon)]$ with $B/\varepsilon\le d$.
\end{proof}
A uniform $o_{\PP}(1/d)$ rate fails: for $Q=\diag(1,1,0,\ldots,0)$ and one unit time pair, the double positive sample and root values are $S=|\det G_{2\times2}|$ and $\sqrt2$, giving supremum error $2|S-\sqrt2|/[d(S+\sqrt2)]$.

Under \eqref{eq:pc} convergence can be arbitrarily slow even at one scale. In the isotropic Brownian model take $X_i=B_i g_i/\sqrt{u_d}$ with independent Gaussian rows and Bernoulli$(u_d)$ switches, $u_d=1-p_d$, $p_d\downarrow0$, $dp_d/(1+\log d)^2\to\infty$. These rows are isotropic with $\Delta_d(a)\le2p_d+\sqrt{2/(ad)}$. For $h(t)=bt/(1+bt)$ and $\widehat u_d(t)=\int w_t\,\dd\mu_d^{I,X}$,
\[
 \sup_{t>0}|(\Theta_d(t)-\widehat u_d(t))/p_d-h(t)^2|\to0\qquad\text{almost surely}.
\]
The active fraction $m/d=u_d+o(p_d)$ almost surely by binomial concentration. The active block is Gaussian of size $m$, scaled by $1/u_d$; \eqref{eq:two-spectrum-tracking} and \eqref{reb:brownian-quantile} give $\widehat u_d(t)=(m/d)bt/(u_d+bt)+o(p_d)$ and $\Theta_d=h+O(1/d)$ uniformly, with summable failure probabilities on $m\ge d/2$.

\begin{proof}[Proof of \eqref{eq:inverse-modulus}]
Put $R=(T/s)^2$, write $A_i(t)=x_i(t)-c_i(t,t)$ temporarily, and recall $1\le r_i\le R$. The anchor denominator is at least $x_i(s)-c_i(T,s)\ge\delta/(R-1)$ by \eqref{eq:relative-parameter}. Also $A_i(t)\ge A_i(s)/R\ge\delta/[R(R-1)]$ by \eqref{eq:spatial-derivative}. For diagonal data set $f_i=1/A_i$. The lower bound gives $\|f_1-f_2\|_\infty=O(e)$ and $\operatorname{Var}_{[s,T]} f_i=O(1)$ since both $x_i$ and $c_i(t,t)$ increase in $[0,1]$. Integration by parts in $\int_s^t f_i\,\dd x_i$ gives an $O(e)$ difference uniformly in $t$. For section data use the quotient directly. This gives the first bound and $\|y_1-y_2\|_\infty\le Ce$.

The arguments $r_i,y_i$ stay in fixed compact positive intervals since $x_i(T)\le R x_i(s)$ and $y_i(s)=sx_i(s)$. The anchor equations bound the source masses, hence $H_i,G_i$ are uniformly Lipschitz there. As $x_i$ changes by at least $\delta$, each argument traverses an interval of length bounded below. The clock estimate and \eqref{eq:two-copy-transforms} therefore compare both compact transforms above to accuracy $Ce$, each on a real interval of uniformly positive length in a compact subset of $(-1,\infty)$.

For bounded finite measures let $f=F_{\nu_1}-F_{\nu_2}$, with $|f|\le\epsilon$ on such an interval. Every admissible interval contains one of finitely many fixed segments $J$ disjoint from a small disk $|z|\le r$ about zero. Take a bounded rectangle $U\Subset\{\operatorname{Re}z>-1\}$ containing all segments and the disk; $f$ is uniformly bounded there. The two-constants theorem on $U\setminus J$ gives $|f|\le C\epsilon^\gamma$ on the disk, since the harmonic measure of $J$ has positive minimum there. Cauchy's inequality gives $|\Delta m_k|\le Cr^{-k}\epsilon^\gamma$ for all moments. Kong and Valiant~\cite[Proposition~1]{KV} bound $W_1$ for probabilities by $C/n+C^n\max_{1\le k\le n}|\Delta m_k|$. Fill finite measures to a common mass at zero and normalize; undoing the fill costs $|\Delta m_0|$ in bounded--Lipschitz distance. Taking $n$ a small multiple of $\log(1/\epsilon)$ proves \eqref{eq:inverse-modulus}.
\end{proof}
\subsection{Finite-dimensional matching and spectral trimming}\label{reb:matching-section}
Write $u_r=\lambda_r^2$ and order both spectra decreasingly.
\begin{theorem}[The cost of trimming both spectra]\label{reb:matching}
For $t>0$ and integers $0\le m\le d$, set
\[
 C_d(t)=\min_{0\le h\le d,\ \ell\ge0}
 \left\{h+\ell+2t^2e_2(q_{>h})\sum_{r>\ell}u_r\right\},\qquad
 B_{d,m}=\sum_{i>m}q_i^2\sum_{r>m}u_r,
\]
with integer deletion counts $h,\ell$. Then
\begin{equation}\label{reb:matching-cost}
 C_d(t)/4\le d\Theta_d(t)/2\le C_d(t)
\end{equation}
and
\begin{equation}\label{reb:matching-bounds}
 \frac md\frac{4t^2B_{d,m}}{1+4t^2B_{d,m}}
 \le \Theta_d(t)\le \frac{4m}{d}+2t^2B_{d,m}.
\end{equation}
\end{theorem}

\begin{proof}
Assign weight $w_e=2t^2q_iq_ju_r$ to $e=(\{i,j\},r)$, $i<j$. A matching uses each label at most once; fixed sets of $2k$ spatial and $k$ temporal labels admit $(2k)!/2^k$ matchings. Thus $P_d(t)$ normalizes their product weights. With $N$ occupied edges, insertion and deletion give
\begin{equation}\label{reb:insertion}
 \PP_t(e\text{ occupied})=w_e\PP_t(e\text{ free}),\qquad
 d\Theta_d(t)/2=\E_tN=\E_t\sum_{e\text{ free}}w_e .
\end{equation}
For a label set $S$ let $W(S)$ be the weight of edges disjoint from $S$. Deleting the largest weights minimizes $W$ at fixed counts, so $C_d=\min_S(|S|+W(S))$. At most $|S|$ occupied edges meet $S$ and the others have expected count at most $W(S)$. Conversely the labels of a random matching have size $3N$ and $\E_t W(S)=\E_tN$, proving \eqref{reb:matching-cost}, also for countably many pairs by nonnegative summation. At fixed $h$ the temporal minimization caps each $2t^2e_2(q_{>h})u_r$ at one.

Taking $h=\ell=m$ and $2e_2(q_{>m})\le d\sum_{i>m}q_i^2$ gives the upper bound in \eqref{reb:matching-bounds}. If $2N<m$, the free spatial pair weight is at least $e_2(q_{>2N})\ge(m-2N)\sum_{j>m}q_j^2$ by pairing $2N<i\le m<j$, and the free temporal weight is at least $\sum_{r>m}u_r$. Hence $\E_tN\ge2t^2 B_{d,m}\E_t(m-2N)_+$ by \eqref{reb:insertion}, giving the lower bound.
\end{proof}

\section{Endpoint detection under directional small-ball bounds}\label{r22:tightness-section}
Throughout this section $Q_d=I$, $\mu_d=\mu_d^{I,X}$ and the centered isotropic rows are independent. For endpoint detection one may replace the projection condition by
\begin{equation}\tag{SB}\label{r24:sb}
 a=\inf_{d,i,\|h\|=1}\E\min\{|X_{i,d}(h)|^2,1\}>0 .
\end{equation}
Equivalently $\PP(|X_{i,d}(h)|\ge u)\ge v$ uniformly for some $u,v>0$: split the clipped expectation at $a/2$, and conversely use $v\min(u^2,1)$. The quadratic-form variance condition \eqref{eq:qc} implies this with $a\ge1/[8(1+\kappa)]$ by rank-one Paley--Zygmund. For lower covariance estimates on isotropic observations see Yaskov~\cite{YLower}.

\begin{lemma}[Detection through a positive covariance]\label{r22:bulk-detector}
Suppose the $d$ rows satisfy \eqref{r24:sb}. Let $0\preceq C\preceq K^*K$, $L=\tr C>0$, $\tr C^2\le L^2/r$ and $2\le r\le d$ an integer. With $b=a\sqrt{dL}/4$, $\theta=a^4r/(192d)$,
\begin{align}
 \E\mu_d((b,\infty))&\ge\theta,\label{r22:bulk-mean}\\
 \PP\{\mu_d((b,\infty))<\theta/2\}&\le e^{-d\theta^2/8}.\label{r22:bulk-prob}
\end{align}
Under \eqref{eq:qc}, if $r\ge r_0:=\lceil512(1+\kappa)\rceil$, one may take $b=\sqrt{dL}/4$ and $\theta=r/(2^{10}d)$.
\end{lemma}
\begin{proof}
For $Y=[FX_i]_{i\le N}$ put $T=F^*F$, $0<l=\tr T<\infty$ and $k=\min(N,\rank F)$. Clip columns at $\sqrt l$ to obtain $Y'$ and $S=Y'Y'^*\preceq YY^*$. Concavity of $\min(u,1)$ in an eigenbasis of $T$, and independence for cross terms, give
\begin{equation}\label{r22:first-gram}
 \E\tr S\ge aNl,\qquad \E\tr S^2\le Nl^2+N(N-1)\tr T^2 .
\end{equation}
If $k\tr T^2\le l^2$, the latter is at most $2(Nl)^2/k$. Splitting $\tr S$ at $aNl/(2k)$ and using Cauchy--Schwarz above the split gives
\begin{equation}\label{r22:second-gram}
 \E N_{>aNl/(2k)}(YY^*)\ge\beta k,\qquad \beta=a^2/8.
\end{equation}

Take $m=\lfloor r/2\rfloor$, $n=d-m$, $V=[KX_i]_{i\le m}$. Since $C\preceq K^*K$,
\begin{equation}\label{r22:covariance-domination}
 B=V^*V\succeq A=[\langle C^{1/2}X_i,C^{1/2}X_j\rangle]_{i,j\le m}.
\end{equation}
Applying \eqref{r22:second-gram} with $F=C^{1/2}$ gives $\E\ell\ge\beta m$ for $\ell=N_{>aL/2}(B)$. Choose a measurable orthonormal frame $E$ for this subspace, $B_E=E^*BE$ and $U=VE B_E^{-1/2}$. Conditionally, \eqref{r22:second-gram} with $F=U^*$, $N=n$, $l=k=\ell$ gives $\E_2 N_{>an/2}(W^*W)\ge\beta\ell$ for $W=[U^*X_j]_{j>m}^*$ (zero count if $\ell=0$). The compressed off-diagonal block $WB_E^{1/2}$ therefore detects at least $\beta^2m\ge a^4r/192$ singular values above $a\sqrt{nL}/2>b$ in expectation. A row edit costs $2/d$ in tail fraction, so bounded differences give \eqref{r22:bulk-prob}.

Under \eqref{eq:qc} take $m=\lfloor r/64\rfloor$ and $\ell=N_{\ge L/2}(B)$ with its corresponding frame. Independence, \eqref{eq:qc} and $\E_2\|U^*X_j\|^4\le\ell^2+\kappa\ell$ give
\[
 \E\|A-LI_m\|_F^2\le m(m-1+\kappa)\tr C^2,\qquad
 \E_2\|W^*W/n-I_\ell\|_F^2\le\ell(\ell+\kappa-1)/n .
\]
Each discarded eigenvalue costs at least one quarter of the squared target. As $4(m+\kappa)\le\min(r,n)/2$, at least half survive in each expectation. Thus the expected count above $\sqrt{dL}/4$ is at least $m/4\ge r/1024$. Cylindrical evaluations follow by finite projection in $L^2$ after Hilbert--Schmidt maps and conditionally along the independent frame, in $L^4$ under \eqref{eq:qc} by rank-one testing.
\end{proof}

\begin{proposition}[Two-sided response detection]\label{r23:response-detection}\label{r23:mass-comp-low}
For independent isotropic rows,
\begin{equation}\label{r23:response-upper}
 \E\mu_d((R,\infty))\le\min\{1,2G_d(1/R)\}\qquad(R>0).
\end{equation}
Under \eqref{r24:sb}, let $g=G_d(t)$, $h=\min(g,1)$, $dh\ge2$, $c_*=a^4/384$ and $b_t=a\sqrt g/(4t)$. Then
\begin{align}
 \E\mu_d((b_t,\infty))&\ge c_*h,\label{r23:response-lower}\\
 \PP\{\mu_d((b_t,\infty))<c_*h/2\}&\le e^{-c_*^2dh^2/8}.\label{r23:response-prob}
\end{align}
Under \eqref{eq:qc}, for $dh\ge4r_0$ one may take $b_t=\sqrt g/(4t)$ and $c_*=2^{-12}$.
\end{proposition}
\begin{proof}
Split at $x_{k,d}=R$. Removing the upper part costs two ranks per pair. For the lower part $T$, independent isotropy gives $\E\|XTX^{\mathsf T}\|_F^2=d(d-1)\|T\|_{\HS}^2$, so the mean tail is bounded by $(2/d)\sum_k\min(1,x_{k,d}^2/R^2)\le2G_d(1/R)$. For the lower bound use
\begin{equation}\label{r23:regularized-covariance}
 C_t=K_d^*K_d(I+t^2d^2K_d^*K_d)^{-1},\qquad
 L_t=\tr C_t=g/(t^2d).
\end{equation}
Since $\|C_t\|_{\op}\le(t^2d^2)^{-1}$,
\begin{equation}\label{r23:regularized-identities}
 L_t^2/\tr C_t^2\ge L_t/\|C_t\|_{\op}\ge dg .
\end{equation}
Apply Lemma~\ref{r22:bulk-detector} with $r=\lfloor dh\rfloor\ge dh/2$, or under \eqref{eq:qc} with $r=\lfloor dh/2\rfloor\ge dh/4$.
\end{proof}

\begin{lemma}[Trimmed energy and the source response]\label{r22:source-criteria}
For $\varepsilon>0$ put
\begin{equation}\label{r22:trimmed-energy}
 \mathcal E_d(\varepsilon)=d\sum_{k>\lfloor\varepsilon d\rfloor}\lambda_{k,d}^2 .
\end{equation}
Bounded $\limsup_d\mathcal E_d(\varepsilon)$ for every $\varepsilon>0$ is equivalent to $\lim_{t\downarrow0}\limsup_d G_d(t)=0$. Also
\begin{equation}\label{r22:source-zero}
 \eta_d\Rightarrow0\quad\Longleftrightarrow\quad G_d(1)\to0
 \quad\Longleftrightarrow\quad \mathcal E_d(\varepsilon)\to0\text{ for every }\varepsilon>0 .
\end{equation}
\end{lemma}
\begin{proof}
Write $E_m=d\sum_{k>m}\lambda_{k,d}^2$. Cutting the decreasing sequence at $td\lambda_k=1$ gives
\begin{equation}\label{r22:response-energy-bound}
 A_d(t):=\min_{m\ge0}\{2m/d+2t^2 E_m\}
 =\frac2d\sum_k\min\{1,(td\lambda_k)^2\},\qquad G_d(t)\le A_d(t)\le2G_d(t).
\end{equation}
Here the comparison is pointwise in each summand. Testing $m=\lfloor\varepsilon d\rfloor$ gives $G_d(t)\le2\varepsilon+2t^2\mathcal E_d(\varepsilon)$. Conversely a minimizing $m$ satisfies $m/d\le G_d(t)$, hence
\begin{equation}\label{r22:energy-source-bound}
 G_d(t)\le\varepsilon\quad\Longrightarrow\quad t^2\mathcal E_d(\varepsilon)\le G_d(t).
\end{equation}
Use a small fixed $t$ for boundedness and $t=1$ for convergence to zero. Here $G_d(1)=\eta_d(I_*)$.
\end{proof}

\begin{corollary}[The isotropic endpoints]\label{r22:tightness}\label{r22:zero}
Under \eqref{r24:sb}, tightness of the mean laws, almost-sure tightness of the sample laws and bounded $\limsup_d\mathcal E_d(\varepsilon)$ for every $\varepsilon>0$ are equivalent. Moreover
\begin{equation}\label{r22:nullity-main}
 \mu_d\Rightarrow\delta_0\text{ in probability}\quad\Longleftrightarrow\quad
 \mathcal E_d(\varepsilon)\to0\text{ for every }\varepsilon>0 .
\end{equation}
Convergence is then almost sure. Equivalently there is a deterministic skew decomposition
\begin{equation}\label{r22:null-decomposition}
 K_d=R_d+S_d,\qquad\rank R_d=o(d),\qquad d\|S_d\|_{\HS}^2\to0 .
\end{equation}
These equivalences also hold under \eqref{eq:pc}.
\end{corollary}
\begin{proof}
Under \eqref{r24:sb}, the preceding lemma and \eqref{r23:response-upper} give mean tightness in the bounded case and convergence to zero in the vanishing case; use $G_d(t)\le\max(1,t^2)G_d(1)$. These conclusions hold almost surely by the row rank bound
\begin{equation}\label{r22:tail-concentration}
 \PP\{|\mu_d((R,\infty))-\E\mu_d((R,\infty))|>u\}\le2e^{-du^2/2}
\end{equation}
at rational $R,u>0$. Almost-sure tightness implies mean tightness by bounded convergence on the supremum of the tails. If the mean tightness condition on $G_d$ fails, choose $t_j\downarrow0$ and $d_j\uparrow\infty$ with $G_{d_j}(t_j)\ge\varepsilon>0$. By \eqref{r23:response-lower}--\eqref{r23:response-prob} a fixed fraction $\delta>0$ eventually lies above $R_j=a\sqrt\varepsilon/(4t_j)\to\infty$, almost surely. In particular on this event the normalized log determinants satisfy
\begin{equation}\label{r22:escaping-log}
 d_j^{-1}\log\det(I+zXK_{d_j}X^{\mathsf T})\ge\tfrac\delta2\log(1+z^2R_j^2)\to\infty\quad(z>0).
\end{equation}
If $G_d(1)\ge\varepsilon$ on a subsequence, the same lower mean bound at $t=1$ rules out convergence in probability to $\delta_0$. Under \eqref{eq:pc}, Theorem~\ref{thm:main} and the two-spectrum criterion with $Q=I$ give the equivalences directly. Finally truncate the first $\lfloor\varepsilon_d d\rfloor$ skew pairs with $\varepsilon_d\downarrow0$ sufficiently slowly to obtain \eqref{r22:null-decomposition}. Conversely best rank approximation gives $2\mathcal E_d(\varepsilon)\le d\|S_d\|_{\HS}^2$ eventually for each $\varepsilon>0$.
\end{proof}

Fair independent multipliers $R_i\in\{0,\sqrt2\}$ on Gaussian rows preserve isotropy, directional small-ball bounds and the mean polynomial. For $j\lambda_j\to c>0$, the active Gaussian block is doubled and its binomial size is asymptotic to $d/2$. Proposition~\ref{prop:rate} gives summable conditional tails, hence almost-sure limit $\tfrac12\delta_0+\tfrac12\mathcal L(2\pi c|C|)$ for standard Cauchy $C$. To see the order of tests in \eqref{eq:pc-stat}, take $(I-E_i)g_i+R_iE_i g_i$ with $E_i$ of equal even rank $\sim d/8$. On rank-$p$ $P$ the conditional squared norm has mean $p\pm\tr(PE_i)$ and variance at most $8p$. Orthogonal $E_i$ give average relative deviation at most $\sqrt{8/p}+1/d$, while coincident $E_i=P$ give at least one by conditional Jensen.

\section{Dilation, fixed kernels and ordered spectra}\label{r21:sec-dilation}
Keep the isotropic notation of Appendix~\ref{r21:canonical-section}. For distinct positive $r_j\to1$, consider
\begin{equation}\tag{DC}\label{r21:dc}
 G_{\lfloor r_jd\rfloor}(t)-r_j^{-1}G_d(r_jt)\longrightarrow0
 \qquad(j\ge1,\ t>0).
\end{equation}

\begin{theorem}[One observation under dilation]\label{r21:coherent}
Under \eqref{eq:pc} and \eqref{r21:dc}, the following are equivalent for fixed $z_0>0$: convergence of $\mu_d$ in probability to a deterministic finite-line law, convergence $G_d(t)\to bt$ for all $t>0$ and some $b\ge0$, and convergence $\widehat q_d(z_0)\to q_0>0$ in probability to a constant. Then $b=(q_0^{-1}-1)/z_0$, $\eta_d\Rightarrow 2b\,\dd x/[\pi(1+x^2)]$, and $\mu_d\Rightarrow\HC_b$ almost surely. Under \eqref{r21:li} one further equivalent test is $F_d(z_0)\to L<\infty$ in probability to a constant, giving $b=(e^L-1)/z_0$. The limits $\widehat q_d(z)\to(1+bz)^{-1}$ and, under \eqref{r21:li}, $F_d(z)\to\log(1+bz)$ hold in $L^1$ and almost surely locally uniformly in $\log z$.
\end{theorem}
\begin{proof}
A finite-line limit gives $q_0>0$. Conversely $q_d(z_0)\to q_0>0$ bounds the source masses by $G_d(z_0q_d)=1-q_d$ and $g_y\ge\min(1,y^2)$. Local equicontinuity thus gives $G_d(u)\to1-q_0$ at $u=z_0q_0$. Every compactified cluster transform $G$ satisfies $G(r_ju)=r_j(1-q_0)$ by \eqref{r21:dc}, hence $G(t)=bt$ by analyticity. This transform gives the stated source and Cauchy law by \eqref{r21:weyl-source} and Corollary~\ref{r21:moving}; \eqref{r21:moving-equivalent} gives the claimed response convergence.
If $F_d(z_0)\to L$ in probability then the mean laws are tight, since
\begin{equation}\label{r22:log-forces-tightness}
 \mu_d((R,\infty))\le\min\{1,2F_d(z_0)/\log(1+z_0^2R^2)\}.
\end{equation}
Corollary~\ref{r21:moving} gives source tightness. With \eqref{r21:li}, Theorem~\ref{r21:logarithms} gives $V_d(z_0)\to L$. For any source cluster point $\eta$, its dilate $\eta^{[r]}=(D_r)_*[r(1+s^2)(1+r^2s^2)^{-1}\eta(\dd s)]$ has response and potential $G_\eta(rt)/r,\Phi_\eta(rt)/r$. By \eqref{r21:dc} and injectivity, $h(r)=V_{\eta^{[r]}}(z_0)$ equals $L$ at every $r_j$. Substituting $w=rq$ gives
\begin{equation}\label{r21:log-envelope}
 rh(r)=\min_{w>0}\{w-r\log(w/r)-r+\Phi_\eta(z_0w)\},\quad
 w_r+G_\eta(z_0w_r)=r .
\end{equation}
The minimizer $w_r$ is analytic by the inverse function theorem, and $(rh)'=\log(r/w_r)$. Thus $h\equiv L$, $w_r=re^{-L}$, and $G_\eta(t)=(e^L-1)t/z_0$. Conversely $G_\eta(t)=bt$ gives $V_\eta(z)=\log(1+bz)$, so Theorem~\ref{r21:logarithms} gives the logarithmic limits.
\end{proof}

\begin{lemma}[Detecting the total source mass]\label{r23:mass-coercivity}
Under \eqref{eq:pc}, if $m_d=G_d(1)\to\infty$, then $\widehat q_d(z)\to0$ almost surely and in $L^1$ for every $z>0$. In particular $\liminf_d\E\widehat q_d(z_0)>0$ implies $\sup_dm_d<\infty$. Under \eqref{eq:qc}, for $m_d>0$, $0<\delta\le1$ and $d\ge16/\delta$,
\begin{equation}\label{qc:mass-bound}
 \E\widehat q_d(z)\le\delta+\frac{32(1+\kappa/d)}{\delta m_d}
 +\frac{64\kappa}{\delta d}+(1+z^2\delta^2m_d/256)^{-1}.
\end{equation}
\end{lemma}
\begin{proof}
The qualitative claims follow from $1-q_d(z)\ge m_d\min(1,z^2q_d(z)^2)$ and Theorem~\ref{thm:main}. For the bound, put $S=K_d^*K_d$, $C=S(I+d^2S)^{-1}$, $L=\tr C=m_d/d$, so $C\le S$, $\tr C^2\le L/d^2$. For the $i$th principal minor $B_i$ of $A_d^X$, put $D_i=(I+z^2B_i^*B_i)^{-1}$, $V_i=(K_dX_j)_{j\ne i}$ and $M_i=z^2V_iD_iV_i^*$. The Schur complement gives $R_{ii}=[(I+zA_d^X)^{-1}]_{ii}=(1+\langle X_i,M_iX_i\rangle)^{-1}$, and interlacing gives $|\tr D_i-d\widehat q_d|\le3$. For the Gram matrix $A_i$ of $(C^{1/2}X_j)_{j\ne i}$, independence and \eqref{eq:qc} give
$\E\|A_i-LI\|_F^2\le(d-1)(d-2+\kappa)L/d^2$. Except on an event of probability at most $32(1+\kappa/d)/(\delta m_d)$, at most $\delta d/8$ eigenvalues of $A_i$ lie below $L/2$. If $\tr D_i\ge\delta d/2$, at least $\delta d/4$ eigenvalues of $D_i$ are at least $\delta/4$. On the good event $V_i^*V_i\ge A_i$; intersecting the latter large-eigenvalue space of $D_i$ with the inverse image under $D_i^{1/2}$ of the large-eigenvalue space of $A_i$ shows $M_i\ge(z^2L\delta/8)P$ for a projection $P$ of rank $k=\lfloor\delta d/8\rfloor\ge\delta d/16$, determined by the other rows. Conditional \eqref{eq:qc} gives $\PP_i(\|PX_i\|^2<k/2)\le4\kappa/k$. Finally $\widehat q_d\ge\delta$ implies $\tr D_i\ge\delta d/2$. Average $R_{ii}$ using the preceding estimates on this latter event; the contribution from $\widehat q_d<\delta$ is at most $\delta$.
\end{proof}

\begin{corollary}[Fixed-kernel endpoints and logarithmic observation]\label{r22:fixed-boundaries}
For fixed $K$, under PC or SB, tightness of the means and almost-sure sample tightness are each equivalent to $A=\sup_k k\lambda_k<\infty$, i.e.\ $K\in\mathcal S_{1,\infty}$. Collapse to $\delta_0$ in probability or almost surely is equivalent to $k\lambda_k\to0$, i.e.\ $K\in\mathcal S_{1,\infty}^0$. Boundedness in probability of $F_d(z_0)$ at one $z_0>0$ forces $A<\infty$. Under PC, $F_d(z_0)\to L<\infty$ in probability to a constant is equivalent to the conditions of Corollary~\ref{thm:fixed}; then Theorem~\ref{r21:coherent} applies with $b=\pi c=(e^L-1)/z_0$.
\end{corollary}
\begin{proof}[Proof of both fixed-kernel corollaries]
The bounds $n^2\lambda_{2n}^2\le n\sum_{k>n}\lambda_k^2\le(\sup_{k>n}k\lambda_k)^2$ give both endpoint criteria by Corollary~\ref{r22:tightness}; bounded logarithms give tightness by \eqref{r22:log-forces-tightness}. Fixed $K$ satisfies \eqref{r21:dc}: $G_n(t)=(d/n)G_d(nt/d)$ and $0\le uG'_d(u)\le2G_d(u)$ bound the rounding error by $O_r(G_d(rt)/d)=o(1)$, since $G_d(u)/d=2\sum_j u^2\lambda_j^2/(1+d^2u^2\lambda_j^2)\to0$ by domination.

A positive limiting complementary response gives tightness by Theorem~\ref{r21:coherent}, hence $A<\infty$. Lemma~\ref{r21:weak-profile} identifies the source in \eqref{r21:weyl-source} with the profile limit $d\lambda_{\lceil dx\rceil}\to c/x$, $c=b/\pi$, pointwise by monotonicity. This is equivalent to $j\lambda_j\to c$ by rescaling. Finally $A<\infty$ implies \eqref{r21:li} by Proposition~\ref{r21:weak-rates}, giving the logarithmic criterion.
\end{proof}
For fixed $K$, \eqref{eq:free-time-law} gives $\varphi_d(w)=\varphi_1(w/d)$, hence $\pi_d=(D_d)_*(\pi_1^{\boxplus1/d})$. The free L\'evy measure of $\pi_1$ is $\sum_{\lambda_k>0}(\delta_{\lambda_k}+\delta_{-\lambda_k})$ \cite{BP,BG,AH}.

For $\lambda_k=k^{-p}$, $1/2<p<1$, $G_d(1)\ge\lfloor d^{1/p}\rfloor/d\to\infty$, illustrating the condition $q_0>0$. Without \eqref{r21:dc}, $\lfloor\gamma d\rfloor$ pairs of size $B/d$ ($0<\gamma<1/2$) give zero atom $1-2\gamma$ and logarithmic limit $\min_{q>0}\{q-\log q-1+\gamma\log(1+B^2z_0^2q^2)\}$ by Theorem~\ref{r21:logarithms}. Since $1-2\gamma\le q\le1$ at the minimum, this limit increases continuously and strictly from zero to infinity with $B$. Alternating two choices with distinct $\gamma$ and the same positive minimum preserves only the logarithmic observation.

\subsection{Ordered spectral recovery under domination}\label{app:ordered}
For nonnegative decreasing profiles with $\int_0^\infty\log(1+a^2)\,\dd x<\infty$, identified almost everywhere, put
\begin{equation}\label{r21:profile-source}
 \eta_a(B)=2\int_{\{a>0\}}\mathbf1_B(a(x))\frac{a(x)^2}{1+a(x)^2}\,\dd x,
 \quad \rho_*(a,b)=\int_0^\infty e^{-x}\min\{1,|a(x)-b(x)|\}\,\dd x .
\end{equation}
Then $a_d(x)=d\lambda_{\lceil dx\rceil,d}$ gives $\eta_{a_d}=\eta_d$. Use subscripts $a$ for the corresponding response, complementary root, potential minimum, and output law. Write $\Phi_a(t)=\int\log(1+t^2a^2)\,\dd x$.

\begin{lemma}[Compact weak-trace profiles]\label{r21:weak-profile}
For $A<\infty$ the class $0\le a(x)\le A/x$ of decreasing profiles is compact in $\rho_*$. The map $a\mapsto\eta_a$ is a homeomorphism onto its image. The source masses are uniformly bounded and tight, the logarithmic tails in \eqref{r21:li} vanish uniformly, and
\begin{equation}\label{r21:profile-bounds}
 \max\{\Phi_a(t),G_a(t)\}\le\pi At,\quad
 \eta_a((M,\infty))\le2A/M,\quad \eta_a([0,\eps])\le4A\eps .
\end{equation}
The Weyl profile has
\begin{equation}\label{r21:weyl-source}
 \eta_{c/x}(\dd s)=\frac{2c}{1+s^2}\,\dd s,\quad G_{c/x}(t)=\pi ct,\quad \mu_{c/x}=\HC_{\pi c}.
\end{equation}
\end{lemma}
\begin{proof}
Helly selection gives compactness. Source continuity follows by domination by $2\min(1,A^2/x^2)$ in \eqref{r21:profile-source}, and the lengths $|\{a>y\}|=\frac12\int_{(y,\infty)}(1+s^{-2})\,\dd\eta_a(s)$ determine $a$. Integrating the envelope $A/x$ gives the first bound in \eqref{r21:profile-bounds} and \eqref{r21:li} as in Proposition~\ref{r21:weak-rates}. The tail bounds follow from $|\{a>M\}|\le A/M$ and splitting $2\int_{a\le\eps}a^2\,\dd x$ at $A/\eps$. Substitution $s=c/x$ gives \eqref{r21:weyl-source}.
\end{proof}

\begin{proposition}[One-response recovery under spectral domination]\label{r21:ordered}
Suppose $0\le a\le b\le A/x$ are decreasing profiles and $A,z_0>0$. Put $\tau=z_0/(1+\pi Az_0)$, $D(t)=\Phi_b(t)-\Phi_a(t)$ and $\Delta=V_b(z_0)-V_a(z_0)$. Then
\begin{align}
 0\le D(\tau)&\le\Delta\le D(z_0),\label{r21:ordered-source}\\
 \rho_*(a,b)^2&\le(4A+\tau^{-2})\Delta.\label{r21:ordered-inverse}
\end{align}
For each fixed $A,z_0$ the exponent $1/2$ for $\rho_*$ is optimal uniformly over ordered pairs.
\end{proposition}
\begin{proof}
The root equation and \eqref{r21:profile-bounds} give $z_0 q_b(z_0)\ge\tau$. Testing each objective at the other's minimizer bounds $\Delta$ between $D(z_0 q_b(z_0))$ and $D(z_0 q_a(z_0))$. Since $t^2/(1+t^2v)$ increases in $t$, its integral over $a^2\le v\le b^2$ does too, proving \eqref{r21:ordered-source}.

Put $h=\min(1,b-a)$, $E=\{b\ge1/\tau\}$ and $\ell=\log[(1+\tau^2b^2)/(1+\tau^2a^2)]$. Integrating $\tau^2/(1+\tau^2v)$ gives $\ell\ge (b-a)/(2b)\ge xh/(2A)$ on $E$ and $\ell\ge\tau^2 h^2/2$ on $E^c$. For $0\le f\le1$, Fubini gives $(\int f)^2=2\int f(x)\int_0^x f(y)\,\dd y\,\dd x\le2\int xf(x)\,\dd x$. Apply this to $f=e^{-x}h\mathbf1_E$, and Cauchy--Schwarz on $E^c$, where $\int e^{-2x}\,\dd x=1/2$ on the full half-line. The contributions to $\rho_*$ are bounded by $(4A\int_E\ell)^{1/2}$ and $(\tau^{-2}\int_{E^c}\ell)^{1/2}$; Cauchy--Schwarz gives \eqref{r21:ordered-inverse}. For $a=0$, $b=\eps\mathbf1_{(0,1)}$ and $0<\eps<\min(A,1)$, $\rho_*=(1-e^{-1})\eps$ and $\Delta\le D(z_0)=\log(1+z_0^2\eps^2)$, proving optimality.
\end{proof}

For arrays with ordered profiles $a_d\le b_d\le A/x$ under \eqref{eq:pc} and any joint coupling, put $U_d=F_d^Y(z_0)-F_d^X(z_0)$. Then $U_d\to0$ in probability, $\rho_*(a_d,b_d)\to0$, and $d_{\mathrm L}(\mu_d^X,\mu_d^Y)\to0$ in probability are equivalent and imply almost-sure convergence. Indeed, Proposition~\ref{r21:weak-rates} and \eqref{r21:ordered-inverse} give the profile estimate, and compactness (Lemma~\ref{r21:weak-profile}), Corollary~\ref{r23:two-array-recovery} and Theorem~\ref{r21:logarithms} give the converse implications by subsequences. Under \eqref{eq:qc} the same estimate reads
$\rho_*(a_d,b_d)^2\le(4A+\tau^{-2})(\E|U_d|+C_{A,\kappa,z_0}d^{-1/4}\log^{3/4}(d+2))$.
Compression $PKP$ supplies ordered pairs by the singular-value inequality.
\section{Attainable mesh profiles and logarithmic loss}\label{app:mesh}
Let $\Pi_b(\nu)$ have complementary response $q+B_\nu(zq)=1$ at $Q_d=I$. Realize mass $a$ at $0<u<\infty$ by $\sim ad/u$ equal cells occupying time fraction $a$; use $d^2$ cells at zero and one at infinity. Atomic approximation and diagonalization realize every $\nu\in\mathcal P(I_*)$. By compactness and Corollary~\ref{prop:joint-mesh}, $\Pi_b$ is a homeomorphism onto the attainable laws. Under \eqref{eq:pc}, $d_{\mathrm L}(\mu_d^{I,X},\Pi_b(\nu_d))\to0$ almost surely, so the sample and profile limits correspond. Put $L=B_\nu(\infty)=\int u^{-1}\,\dd\nu$ with extended reciprocals. Then $\Pi_b(\nu)(\{0\})=(1-L)_+$ since $q(\infty)\ge1-L$ with equality whenever $q(\infty)>0$ by $zq(z)\to\infty$.
For $\rho$ on $[0,\infty)$ with positive mass off zero, set $a=bt$, $h=H_\rho(a)$, $L=\int(1\wedge u)\,\dd\nu$ and $J=\int hu/(1+hu)\,\dd\nu$ (one at infinity). Then for every $t>0$,
\[
 A_\rho(h)J\le h-y(t)\le2hL,\qquad h-y(t)\le A_\rho(h)J/(1-2L)^2\quad(L<1/2).
\]
Thus $(h-y(t))/J\to A_\rho(h)$ for fixed $t$ as $L\to0$ with $J>0$; the finite root error is $4/d$. For $\nu\ne\delta_\infty$ write $y=H_\rho(as)$, $0<s\le1$. The crossing $\int(s+yu)^{-1}\,\dd\nu=1$ gives $1-s=\int yu/(s+yu)\,\dd\nu$. Concavity gives $h\le y/s\le h/s$, hence $J\le1-s\le J/s$. Since $(s+yu)^{-1}\ge(2-s)(1-u)_+$, $1-s\le2L$. The factor $h-C_\rho(a,as)$ in $h-y=(1-s)[h-C_\rho(a,as)]$ is at most $h$ and between $A_\rho(h)$ and $A_\rho(h)/s$ by $1\le(1+aq)/(1+saq)\le1/s$. For $\nu=\delta_\infty$ take $y=0$, $J=L=1$.

\begin{corollary}[Logarithmic loss of a partition]\label{r21:mesh}
Put $\mathcal L_d=\int(1\wedge u)\,\dd\nu_d$, $\tau=z/(1+bz)$, $m_z=(b\tau)^2/[2(1+b\tau)]$,
$M_z=\max((bz)^2/2,bz)$, $\epsilon_d(t)=d^{-1}\log(2(1+btd))$. For $z>0$,
$m_z\mathcal L_d\le\log(1+bz)-V_d(z)\le M_z\mathcal L_d+\epsilon_d(z)$.
Under \eqref{eq:qc}, $\E|F_d(z)-\log(1+bz)|$ satisfies the same bounds with an additive
error $e_d=C_{T,\kappa,z}d^{-1/4}\log^{3/4}(d+2)$ on either side.
Under \eqref{eq:pc} the cell-size criterion in Corollary~\ref{prop:joint-mesh} is equivalent to each of the following convergences in probability:
$F_d(z)\to\log(1+bz)$ at one $z>0$, the complementary response tending to $(1+bz)^{-1}$ at one $z>0$, and $\mu_d^{I,X}\Rightarrow\HC_b$.
These convergences then hold almost surely, locally uniformly in $\log z$ for the scalar statistics.
\end{corollary}
\begin{proof}
Write $x_j=btdp_j$, $\Psi_d=d^{-1}\sum_j\log(1+x_j)$. The even product is $\prod_j(1+x_j)(1+r)/2$, $r=\prod_j(1-x_j)/(1+x_j)$. If $r<0$, some $x_j>1$ gives $|r|\le(x_j-1)/(x_j+1)$, hence $0\le\Psi_d-\Phi_d\le\epsilon_d(t)$ in all cases. With $a=bt$, $bt-\Psi_d=\int f_a\,\dd\nu_d$ where $f_a(u)=a-u^{-1}\log(1+au)=a^2u\int_0^1 s/(1+aus)\,\dd s$. On $[0,1]$, $a^2u/[2(1+a)]\le f_a(u)\le a^2u/2$; on $[1,\infty]$, $f_a$ increases to $a$. Testing each objective at the other's minimizer bounds the loss between $b\tau-\Phi_d(\tau)$ and $bt-\Phi_d(t)$ at $t=zq_d(z)\le z$. Cell compression gives $\lambda_j(K_\pi)\le T/((2j-1)\pi)$, so Proposition~\ref{r21:weak-rates} supplies the sample bounds. Since $\mathcal L_d\to0\iff\nu_d\Rightarrow\delta_0$, Corollary~\ref{prop:joint-mesh} gives the equivalences; equicontinuity in $\log z$ gives local uniformity.
\end{proof}

For $n_d/d\to a\in(0,\infty)$ equal cells give $\nu=\delta_{1/a}$, $q+bzq/(1+bzq/a)=1$ and zero mass $(1-a)_+$; $a=\infty,0$ give $\HC_b,\delta_0$. Using $\lceil\sqrt d\rceil$ equal cells on one time half and $\lceil d^2\rceil$ on the other gives $\nu=(\delta_0+\delta_\infty)/2$ and $\HC_{T/4}$ despite $n_d/d\to\infty$, $\max h_j\to0$. With $p_1=d^{-1/3}$ and $d^2$ equal remaining cells, $\mathcal L_d\le d^{-1/3}+d^{-1}$ gives bulk recovery although It\^o isometry and cell independence give $d^{-1}\E\|A^d-A^{d,\pi}\|_F^2=(d-1)\sum h_j^2/4\sim T^2d^{1/3}/4$. For the first $m_d$ pairs $\lambda_j=c/(j-1/2)$, $c>0$, $m_d/d\to\gamma$, Riemann sums give $G_d(t)\to2ct\arctan(\gamma/(ct))$. Thus $\HC_{\pi c}$ occurs exactly when $m_d/d\to\infty$, and $m_d/d\to0$ gives $\delta_0$; here $c=T/(2\pi)$ for Brownian area.
\paragraph{Finite samples after exchanging the weights.}
After exchange, \eqref{eq:two-spectrum-tracking} and the common padded polynomial give discrepancy $\le2\beta_D$ outside probability $4D^{-3}$, hence eventually almost surely ($D=\max(d,n_d)$ and $\beta_D$ is the bound there). At $d=N=3$, $Q=I$, write $g\times h=RU$ for independent standard Gaussian vectors. Exchange takes the squared positive singular value $b^2e_2(p)R^2$ to $3b^2R^2(p_2p_3U_1^2+p_1p_3U_2^2+p_1p_2U_3^2)$. The uniform direction $U$ is independent of $R$, so the factor in parentheses has mean $e_2(p)/3$ and positive variance for unequal positive $p_i$, separating fourth moments.
\section*{Declaration of AI Use}
The author used GPT interactively in developing and refining the mathematical arguments and the presentation of this paper. The author takes full responsibility for the correctness and mathematical content of the manuscript.
\begingroup
\renewcommand{\bibliofont}{\normalsize}

\endgroup
\end{document}